\documentclass[10pt]{amsart}
\usepackage{amssymb,amsmath,amsfonts,amscd, mathrsfs}
\usepackage{times,graphicx,epsfig}
\usepackage{soul}

\theoremstyle{plain}
\newtheorem{theorem}{Theorem}[section]
\newtheorem{proposition}[theorem]{Proposition}
\newtheorem{definition}[theorem]{Definition}

\newtheorem{corollary}[theorem]{Corollary}
\newtheorem{remark}[theorem]{Remark}

\def\g{\mathfrak{g}}

\def\e{\mathbf{e}}
\newcommand{\Rn}{\mathbb R^n}

\newcommand{\rn}[1]{{\mathbb R}^{#1}}
\newcommand{\R}{\mathbb R}

\newcommand{\G}{\mathbb G}

\newcommand{\N}{\mathbb N}

\newcommand{\supp}{\mathrm{supp}\;}

\newcommand{\scal}[2]{\langle {#1} , {#2}\rangle}

\newcommand{\Scal}[2]{\langle {#1} \vert {#2}\rangle}

\newcommand{\res}{\mathop{\hbox{\vrule height 7pt width .5pt depth 0pt
\vrule height .5pt width 6pt depth 0pt}}\nolimits}

\newcommand{\ccheck}{{\vphantom i}^{\mathrm v}\!\,}
\newcommand{\mc}{\mathcal }

\newcommand{\divg}{\mathrm{div}_H\,}

\usepackage[usenames]{color}

\begin{document}

%\today

\bigskip

\title[The divergence of  vector fields vanishing at infinity on Carnot groups]
{The distributional divergence of horizontal vector fields vanishing at infinity on Carnot groups}

\author[Annalisa Baldi, Francescopaolo Montefalcone]{
Annalisa Baldi\\
Francescopaolo Montefalcone
}

\begin{abstract}
We define a $BV$-type space  in the setting of Carnot groups (i.e., simply connected Lie groups  with stratified nilpotent Lie algebra) that  allows one to characterize all  distributions $F$ for which there exists a continuous horizontal vector field $\Phi$, vanishing at infinity, that solves the equation $\mathrm{div}_H \Phi=F$.
This generalize to the setting of Carnot groups some results by De Pauw and Pfeffer, \cite{DePauw-Pfeffer}, and  by De Pauw and Torres, \cite{DePauw-Torres}, for the Euclidean setting.
\end{abstract}
 


\maketitle

\tableofcontents

\section{Introduction}\label{introduction}

In their celebrated 2003 paper \cite{BB2003}, Bourgain and Brezis studied a problem concerning the equation $\mathrm{div}\,Y=f$ for $L^p$-periodic functions $f$ defined on $\R^n$. Among their results, they considered the limiting case $p=n$ and proved that there exists a vector field $Y$ solving the equation and that belongs to $L^\infty$. 
To attack the above problem,  they  started by using special vector fields  of the form $Y=\nabla u$, thus considering the problem $\Delta u=f$. This method for $1<p<\infty$ yields a solution $u\in W^{2,p}$ and, consequently, a solution $Y\in  W^{1,p}$. Unfortunately, in the limiting case $p=n$,  the fact that  $Y\in  W^{1,p}$ does not imply directly that $\nabla u$ belongs to $L^\infty$, since $ W^{1,p}$ is not contained in $L^\infty$. Despite this, they proved that in a suitable class of periodic functions on $\R^n$
%It is well known that $Y=(-\Delta)^{-1}\mathrm{curl}\,f$ is a solution of \eqref{system intro}. Then, by the Calder\'on-Zygmund theory we can say that
%$$
%\|\nabla Z\|_{L^{p}(\R^3)}\le C_p\|f\|_{L^p(\R^3)}\,, \quad \mathrm{for}\quad {1<p<\infty}.
%$$
%Then, by  Sobolev inequality, if  $1<p<3$ we have:
%$$
%\|\Z\|_{L^{p*}(\R^3)}\le\|f\|_{L^p(\R^3)}\,,
%$$ where $\frac{1}{p*}=\frac{1}{p}-\frac{1}{3}$. When we turn to the case $p=1$ the first inequality fails. The second remains true. This is exactly the result proved in \cite{BB2004} by Bourgain \& Brezis.
%
there exists indeed a {\it continuous} vector fields $Y$ that solves the equation $\mathrm{div}\,Y=f$ and such that 
\begin{equation}\label{1in}
	\|Y\|_{L^{\infty}}\le \mathscr C(n)\|f\|_{L^n},
\end{equation} where $\mathscr C(n)$ denotes a dimensional constant.
 The continuity of $Y$ is a key point in their  proof, which relies on the Sobolev embedding of both spaces $W^{1,1}$ and $BV$ into $L^{n/n-1}$, and on a duality argument. The proof itself was not constructive. As a matter of fact, the authors showed that there are no bounded linear operators $K$ from the space of $L^n$-periodic functions to $L^\infty$ such that $\mathrm{div} (K f)=f$ in the  distributional sense. 
 Thus,  inequality \eqref{1in}  cannot follow from a representation formula for solutions to the  equation under study.
 After the paper \cite{BB2002} was written, a huge literature appeared concerning  equations such as
 \begin{equation}\label{top}
 \mathrm{div}\,Y=F.
 \end{equation}
 Among them we quote  \cite{DePauw-Pfeffer}, where the authors considered the problem  in a more general framework,   finding necessary and sufficient conditions on $F$ in order to get a continuous weak solution of \eqref{top}. Moreover, they introduced the notions of charge and strong charge, which  originated from their researches on generalized Riemann integrals and Gauss-Green theorems; see  \cite{DePauw-Pfeffer} and references therein.
 
 We remind the reader that a  distribution $F\in\mc D'(\R^n)$ is said 
 a {\it flux} if
 the equation  \eqref{top} has a continuous solution, i.e., if there exists a vector field
 $Y\in C(\R^n;\R^n)$ such that 
$$F(\varphi)=-\int_{\R^n}
\langle Y(x) , \nabla\varphi(x)\rangle\,dx \qquad\forall\, \varphi\in\mc D(\R^n).$$

A linear functional $F: \mc D(\R^n)\longrightarrow\R$
is called
a {\it charge} in $\R^n$ if $\lim_{i\to+\infty}F(\varphi_i ) = 0$ for every sequence
$\{\varphi_i\}_{i\in\N}\subset\mc D(\R^n)$ such that
$$\lim_{i\to+\infty}\|\varphi_i \|_{L^1} = 0\qquad \mbox{and}\qquad \sup_{i}(\|\nabla\varphi_i \|_{L^1} + \|\varphi_i \|_{L^\infty}) <\infty;$$ see Definition 2.3 in \cite{DePauw-Pfeffer}. On the other hand, the linear functional $F: \mc D(\R^n)\longrightarrow\R$ is said 
a {\it strong charge} in $\R^n$  if $\lim_{i\to+\infty} F(\varphi_i ) = 0$ for every   sequence
$\{\varphi_i\}_{i\in\N}\subset \mc D(\R^n)$ such that
$\lim_{i\to+\infty}\|\varphi_i \|_{L^1} = 0$ and $\sup_{i}\|\nabla\varphi_i \|_{L^1} <+\infty$.
The linear spaces of all fluxes,  charges   and strong charges in $\R^n$ are denoted, respectively, by
$\mc F(\R^n)$, $\mathbf{Ch}(\R^n)$, and $\mathbf{Ch}_s(\R^n)$. It is observed in \cite{DePauw-Pfeffer} that, in principle, 
$\mc F(\R^n)\subset \mathbf{Ch}(\R^n)\subset \mathbf{Ch}_s(\R^n)\subset \mc D'(\R^n)$ but in the paper the authors show that $\mc F=\mathbf{Ch}_s$. 

We remark that an example of strong charge is given by any distribution associated with a function $f\in L^n_{{loc}}(\R^n)$: this shows the connection  with the  problem studied by Bourgain and Brezis. 
Later on, De Pauw and Torres, \cite{DePauw-Torres},  characterized all functionals $F$ acting linearly on the subspace of $L^{n/(n-1)}(\Rn)$ of all functions whose distributional gradient is a vector valued measure, under a suitable continuity assumption. The requirement on $F$  is connected with the definition of charge vanishing at infinity (see Definition 3.1 in \cite{DePauw-Torres}). As a corollary of their characterization result, De Pauw and Torres proved that given $f\in L^n(\R^n)$ there exists $Y\in C_0(\R^n,\R^n)$ such that $\mathrm{div}\,Y=f$ in the sense of distribution, where $C_0(\R^n,\R^n)$ denotes the space of all continuous vector fields vanishing at infinity.

Starting from the existence result of De Pauw and Torres and adapting Bourgain and Brezis' proof,  Moonens and Picon proved in  \cite{Moonens-Picon}  that if $f\in L^n(\R^n)$, then there exists  $\widetilde Y\in C_0(\R^n,\R^n)$ solving the equation $\mathrm{div}\,\widetilde Y=f$, {\it and such that} $$\|\widetilde Y\|_{L^{\infty}}\le \mathscr C(n)\|f\|_{L^n},$$ where the constant $\mathscr C(n)$ is a dimensional constant independent of $f$.

In this paper we  study, in the setting of  Carnot
groups (i.e., simply connected Lie groups $\mathbb G$, with stratified nilpotent Lie algebra $\mathfrak g$; see, e.g., \cite{BLU},  \cite{FS}, \cite{stein}), an analogous of the equation \eqref{top},  obtaining also a continuity estimate similar to the one  above. 
Carnot groups  are the simpler examples of sub-Riemannian manifolds and play a deep role in studying, in a sub-Riemannian setting, problems arising from differential geometry, geometric measure theory,
subelliptic differential equations,  optimal control theory, mathematical models in
neurosciences and  robotics. Roughly speaking, a sub-Riemannian structure on a smooth $n$-dimensional manifold $M$
is given by a subbundle $H M$ of the tangent bundle $TM$, which defines a family of admissible directions at
any point of $M$. 
The subbundle $H M$ is called the {\it horizontal} bundle. If we endow each fiber $H_x M$  of $H M$ with a scalar product
$\scal{}{}_x$,
there exists a naturally associated distance $d$ on $M$, called Carnot-Carath\'eodory distance, defined as the infimum of the Riemannian length of all {\it horizontal
curves} (i.e., any curve $\gamma:I\to M$ such that $\gamma'(t)\in H_{\gamma(t)}M$ for a.e. $t\in I$)  joining two given points.

In any Carnot group $\G$, the  horizontal subbundle $H\G$ is generated by left translation of 
the first layer of the stratification
of the Lie algebra $\g$, which can be identified with a linear subspace of the tangent space of the group
at the  identity. Moreover, through the Lie group
exponential map, $\G$ can be identified with the Euclidean space $\R^n$, endowed with a polynomial group law, where $n = {\rm dim} \ \mathfrak g$. Notice that the Hausdorff dimension $Q$ of a Carnot group $\G$ turns out to be strictly greater than its  topological dimension.

Horizontal vector fields in  Carnot groups (i.e., smooth sections of the horizontal subbundle $H\G$) are the natural counterpart of vector fields in Euclidean spaces, and there is a well understood notion of horizontal divergence, later denoted as ${\mathrm{div}}_H$.  This fact makes possible to study an equation of the type \begin{equation}\label{Heq}{\text{div}}_H \Phi=F.
\end{equation} More precisely, in this paper we  study  the  notion of charge vanishing at infinity in the setting of Carnot groups, following the lines of \cite{DePauw-Torres},  in connection with the solvability of the equation \eqref{Heq}. 

Our main result is stated in Theorem \ref{main result}, where we prove that
if 
$F\in\mc D'(\G)$, then there exist continuous horizontal vector fields  vanishing at infinity (see Section 2 for precise definitions) that solve \eqref{Heq} in the distributional sense  if and only if  $F$ is a charge vanishing at infinity. As a corollary, if $F\in L^Q(\G)$ (hence, it turns out that $F$ can be regarded as a charge vanishing at infinity), there is a continuous  solution of \eqref{Heq} vanishing at infinity  that in addition satisfies the inequality
\begin{equation}
\label{moonens-piconintro} \|\Phi\|_{L^\infty}\leq \mathscr C(Q) \|F\|_{L^Q},
\end{equation}where $\mathscr C(Q)$ denotes a geometric constant, which is independent of $F$  (see \eqref{moonens-picon}).

The problem of the existence of an $L^\infty$-solution $\Phi$,  and of an inequality like \eqref{moonens-piconintro},  could be formulated in the more general setting of the Rumin complex of intrinsic differential forms on Carnot groups. In fact, 
  horizontal  vector fields  can be identified with intrinsic differential forms of degree $(n-1)$, so that   an estimate like \eqref{moonens-piconintro} can be seen as the first  link of a chain of analogous inequalities for intrinsic differential forms of any degree. 
A similar result, for Rumin's differential forms of any degree, has been recently obtained in the setting of Heisenberg groups in \cite{BFP5}. Nevertheless,  the formulation  of the problem itself, in terms of differential forms of arbitrary degree in   general Carnot groups, is not straightforward at all due to the lack of homogeneity of the Rumin's exterior differential (for an explanation of this phenomenon, see, e.g., \cite{BFTT} p.6).  Thus, one of the motivations of our paper is to attack this kind of problem in  general Carnot groups for  horizontal vector fields (thought of as identified with intrinsic differential forms of degree $(n-1)$), where the Rumin's exterior differential turns out to be always homogeneous. 

\medskip

The paper is organized as follows. 
Precise definitions and basic properties of Carnot groups are discussed in Section 2, together with the notions of horizontal vector measures and horizontal vector fields vanishing at infinity adapted for this setting; see Section 2.1.    
Then, in Section 2.2, we  collect  several results about $BV$  functions in Carnot groups. In Section 3  we introduce and study another $BV$-like space, denoted by $BV^{Q/Q-1}(\G)$ and  defined as the set of all functions in $L^{Q/Q-1}(\G)$ whose distributional gradient (regarded as a measure)  has finite total variation. In Section 4 we study a closed subspace of the dual space of $BV^{Q/Q-1}(\G)$, denoted by  $\mathbf{Ch}_0(\G)$. In particular,  following the lines of  \cite{DePauw-Torres}, we prove that its dual is isomorphic to  $BV^{Q/Q-1}(\G)$. Section 5 contains our main result (see Theorem \ref{main result}) concerning the equation  $\mathrm{div}_H\Phi=F$ (meant in the distributional sense). In particular, we show that this equation admits as a solution a continuous horizontal vector field $\Phi$ vanishing at infinity if and only if $F\in \mathbf{Ch}_0(\G)$. In addition, as a corollary,  we prove an estimate of the type \eqref{moonens-piconintro}; see Corollary \ref{cor main}.

\section{Notation and preliminary results} \label{preliminar} A {\it Carnot group $\G$ of
step $\kappa$}  is a simply connected
Lie group whose Lie algebra ${\mathfrak{g}}$ is finite dimensional, say of
 dimension $n$, and  admits a {\it step $\kappa$ stratification}, i.e.,
there exist linear subspaces $V_1,...,V_\kappa$ such that
\begin{equation}\label{stratificazione}
{\mathfrak{g}}=V_1\oplus...\oplus V_\kappa,\quad [V_1,V_j]=V_{j+1},\quad
V_\kappa\neq\{0\},\quad V_j=\{0\}{\,\,\textrm{ if }\,\,} i>\kappa,
\end{equation}
where $[V_1,V_j]$ denotes the subspace of ${\mathfrak{g}}$ generated by
all commutators of the form $[X,Y]$, with $X\in V_1$ and $Y\in V_j$\, ($j\geq 1$).

 For any  $j=1,\dots,\kappa$, let
$m_j:=\dim {V_j} $ and $h_j:=m_1+\dots +m_j$, where
$h_0=0$ and, clearly, $h_\kappa=n$.  Now choose a basis $\{\e_1,\dots,\e_n\}$ of
$\mathfrak{g}$ adapted to the stratification, i.e.,  
$$\{\e_{h_{j-1}+1},\dots,\e_{h_j}\}\;\text {is a basis of}\; V_j\;\text{
for any}\; j=1,\dots, \kappa.$$
Let $X=\{X_1,\dots,X_{n}\}$ be the set
of left-invariant vector fields  of $\G$ such that
$X_i(e)=\e_i$ $(i=1,...,n)$, where $e$ denotes the identity of $\G$.  
By the stratification hypothesis \eqref{stratificazione},  
 all left-invariant vector fields of $\G$ are  generated  by iterated Lie brackets of 
the subset
$\{X_1,\dots,X_{m_1}\}$: we will
refer to $X_1,\dots,X_{m_1}$ as the {\it generating  vector fields} of
the group.

 The exponential map is a one to one map from $\mathfrak
g$ onto $\G$. Thus, any $x\in\G$ can be written in a unique way as
$x=\exp(x_1X_1+\dots+x_nX_n)$. Using these {\it exponential
coordinates}, we shall identify $x$ with the $n$-tuple $(x_1,\dots,x_n)\in
\R^n$ and, accordingly, $\G$ with $(\Rn,\cdot)$. The explicit
expression of the group operation ``$\cdot$'' follows from the
Campbell-Baker-Hausdorff formula; see \cite{BLU}.
If $j=1,\dots,\kappa$, then set $
x^j:=(x_{h_{j-1}+1},\dots,x_{h_j}) \in\R^{m_j}$. Thus, we can also identify
$x$ with the $\kappa$-tuple $(x^1,\dots, x^\kappa)\in
\R^{m_1}\times\ldots\times\R^{m_\kappa}=\Rn$.

Recall that there are two important families of group automorphisms: left translations and  group dilations. For
any $x\in\G$, the {\it  left translation by $x$}, say $\tau_x:\G\longrightarrow\G$, is the map given by $$ \G\ni z\longmapsto\tau_x z:=x\cdot z. $$ For any $\lambda >0$, the
{\it dilation} $\delta_\lambda:\G\longrightarrow\G$, is  defined as
\begin{equation}\label{dilatazioni}
\delta_\lambda(x_1,...,x_n)=
(\lambda^{d_1}x_1,...,\lambda^{d_n}x_n),
\end{equation} where $d_i\in\N \ (i=1,...,n)$ denotes  the {\it homogeneity} of
the monomial $x_i$ in
$\G$ (see \cite{FS}, Ch.1, par. C), which is given by
\begin{equation}\label{omogeneita2}
d_i=j \;  \quad\text{whenever}\quad \; h_{j-1}+1\leq i\leq h_{j}\quad (j=1,...,\kappa).
\end{equation}
In particular, note that $1=d_1=...=d_{m_1}<
d_{{m_1}+1}=2\leq...\leq d_n=\kappa.$

\medskip

\noindent The Lie algebra $\mathfrak g$ can always be equipped with a scalar product
$\scal{\cdot}{\cdot}$ for which $\{X_1,\dots,X_n\}$ is an orthonormal basis.
 
As customary, we also fix a smooth homogeneous norm $\|\cdot \|$ in $\G$ (see
\cite{stein}, p. 638)
such that the gauge distance $d(x,y):=\|y^{-1}\cdot x\|$  is a 
left-invariant  distance  on $\G$,  in fact
 equivalent to the ``Carnot-Carath\'eodory distance'' (see \cite{ABB}).  
 We set $$B(x,r):=\{y \in  \G; \; d(x,y)< r\}$$ to denote the open $r$-ball centered at $x\in\G$. 
It is well-known that any Haar measure of a Carnot group $\G$ coincides, up to a constant factor, with the standard Lebesgue measure 
$\mathscr L^n$ on $\mathfrak g\cong\rn n$ (notice that we 
just write $dx$ instead of $d\mathscr L^n(x)$ in the integrals). If $A\subset \G$ is a $\mathscr L^n$-measurable set, we will also set $|A|:=\mathscr L^n(A)$.

The {\it homogeneous dimension} $Q$ of the group $\G$ is the number defined as
 \begin{equation}\label{dim_omo}
 Q:=\sum_{j=1}^{\kappa} j \dim V_j.
 \end{equation}
 Since for any $x\in\G$ and $r>0$  we have
\begin{equation}\label{misura palla}
|B(x,r)|=|B(e,r)| = r^Q|B(e,1)|,
\end{equation} the integer
 $Q$ turns out to be the
 Hausdorff dimension of the metric space $(\G,d)$.

\begin{proposition}
\label{legge di gruppo}
The group product ``$\cdot$'' has the form
\begin{equation} x\cdot y=x+y+\mathcal Q(x,y) \qquad
\mbox{for all } x,y\in\Rn, \label{legge di gruppo1}
\end{equation}
where $\mathcal Q=(\mathcal Q_1,\dots,\mathcal
Q_n):\Rn\times\Rn\longrightarrow\Rn$, and any $\mathcal Q_i$ is a homogeneous
polynomial of degree $d_i$ $(i=1,...,n)$ with respect to the intrinsic
dilations \eqref{dilatazioni}, i.e.,
\begin{equation*}\mathcal Q_i(\delta_\lambda x,\delta_\lambda
y)=\lambda^{d_i}\mathcal Q_i(x,y) \qquad \mbox{for all } x,y\in \G.
\end{equation*}  In addition, for  every  $x,y\in\G$ the following hold: 
\begin{gather}
\label{primostrato} \mathcal Q_1(x,y)=...=\mathcal Q_{m_1}(x,y)=0;\\
\label{legge di gruppo3} \mathcal Q_j(x,0)=\mathcal
Q_j(0,y)=0\quad\text{and}\quad\mathcal Q_j(x,x)=\mc Q_j(x,-x)=0 \quad
\text{for}\;m_1< j\leq n;\\ \mathcal Q_j(x,y)=\mathcal
Q_j(x_1,\dots,x_{h_{i-1}},y_1,\dots,y_{h_{i-1}})
\quad\text{for}\quad   h_{i-1}\leq j\leq
h_i\quad (i>1).\label{legge di gruppo2}
\end{gather}
\end{proposition}

It follows from Proposition \ref{legge di gruppo}  that
 $ \delta_\lambda x\cdot\delta_\lambda y=\delta_\lambda (x\cdot y)$  for every 
$x,y\in\G,$ and that the inverse $x^{-1}$ of any $x=(x_1,\dots,x_n)\in\G$ has the form
 $x^{-1}=(-x_1,\dots,-x_n).$ 

\begin{proposition}[see, e.g., \cite{FSSC}, Proposition 2.2] \label{campi omogenei0}
The  left-invariant vector fields $\{X_1,..., X_n\}$ have polynomial coefficients and are of the
form
\begin{equation}\label{campi omogenei}
X_j(x)=\partial_j+\sum_{i>h_l}^n q_{i,j}(x)\partial_i \qquad
\text{for any} \:\; j=1,\dots,n \; \;\text{and}\; \;j\leq h_l\quad(l=1,...,\kappa),
\end{equation} where
$q_{i,j}(x)=\frac{\partial \mc Q_i}{\partial y_j}(x,y){\big|_{y=0}}$.\\

 In particular, if $h_{l-1}<j\leq h_l$, then $q_{i,j}(x)=q_{i,j}(x_1,...,x_{h_{l-1}})$
and $q_{i,j}(0)=0$.
\end{proposition}

The subbundle $H\G$ of the tangent bundle $T\G$ spanned by the
vector fields $\{X_1,\dots,X_{m_1}\}$ is called the {\it horizontal bundle} and plays a particularly important
role in the theory. The fibers of $H\G$ are explicitly given by$$ H_x\G=\mbox{span
}\{X_1(x),\dots,X_{m_1}(x)\} \qquad \forall\, x\in\G .$$ 
For  simplicity of notation, we will henceforth set  $m:=m_1$.

A subriemannian
structure is defined on $\G$ once one endows each fiber $H_x\G$  of the horizontal bundle $H\G$ with a
scalar product $\scal{\cdot}{\cdot}_{x}$; its associated norm is denoted as $|\cdot|_x$. 
When clear from the context, we will drop the subscript $x$, simply writing  $\scal{\cdot}{\cdot}$ and $|\cdot|$.

 From now on, we shall assume that, at any $x\in\G$, the basis $\{X_1(x), \ldots, X_{m}(x)\}$ is
orthonormal  (under the chosen scalar product). 

Measurable sections of the horizontal bundle $H\G$ are called {\it horizontal
sections} (or {\it horizontal
vector fields} ), and vectors in $H_x\G$ are called {\it horizontal vectors}. 

Given a horizontal vector field\footnote{In other words, if $\pi:T\G\to\G$ is the bundle projection map, then $\pi\circ\Phi$ is the identity map.} $\Phi: \G\to H\G$, and since a horizontal frame has already been fixed, we can write $\Phi$ in terms of its $m$ components $\Phi_i: \G\to\R$  $(i=1,\ldots,m)$ along the horizontal frame  $\{X_1, \ldots X_m\} $,  so that
$$
\Phi=\sum_{j=1}^m \phi_j X_j.
$$In other words, we can always assume that  $\Phi=(\phi_1,...,\phi_m)$.

Now, let $f:\G\longrightarrow \R$ be a smooth function, say $f\in C^\infty(\G)$.  The horizontal gradient of $f$ is the horizontal vector field $D_H f$ defined  by $$\scal{D_Hf(x)}{X}_x=df_x(X),\qquad \forall\,x\in\G,\,\,\forall\,  X\in H_x\G.$$
Clearly, with respect to the  the horizontal frame, we can  write $D_Hf=(X_1f,...,X_{m}f)$.
 
 Moreover, if
 $\Phi=(\phi_1,\dots,\phi_{m})$ is a smooth horizontal vector field, say $\Phi\in C^\infty(\G,H\G)$, its horizontal divergence $\divg
 \Phi$ is, by definition,  
  the real valued function
 \begin{equation}\label{divergence def}
 \divg \Phi :=\sum_{j=1}^{m}X_j\phi_j.
 \end{equation}
 
   The same symbols $D_H$ and $\divg$ will be adopted later, when working with the {\it weak} horizontal gradient and divergence operators (intended in the sense of distributions).   
  
 Recall that if $\Omega\subseteq\G$ is an open set, the space of continuous linear functionals on $C^\infty(\Omega)\;(=:\mc E(\Omega))$ is  denoted by $\mc E'(\Omega)$ and the space of
   continuous linear functionals on $C_c^\infty(\Omega)\;(=:\mc D(\Omega))$ is denoted by $\mc D'(\Omega)$. Throughout the paper, we will use the notation $\Scal{\cdot}{\cdot}$ for the duality 
    between $\mc D'(\Omega)$ and $\mc D(\Omega)$ and  also for the duality between $\mc E'(\Omega)$ and $\mc E(\Omega)$ (more generally, the same notation will be used for the duality between other function spaces defined below).

 If $f:\G\longrightarrow\R$, we denote
    by $\ccheck f$ the function given by $\ccheck f(x):=
    f(x^{-1})$. Furthermore, if $T\in\mc D'(\G)$, then $\ccheck T$ will
   denote the distribution defined by $\Scal{\ccheck T}{\varphi}
    :=\Scal{T}{\ccheck\varphi}$ for any test function $\varphi\in\mc D(\G)$.
    
    As in \cite{FS},
we  adopt the following multi-index notation for higher-order derivatives. If $I =
(i_1,\dots,i_{n})$ is a multi--index, we set  
$X^I=X_1^{i_1}\ldots
X_{n}^{i_{n}}$. By the Poincar\'e-Birkhoff-Witt theorem
(see, e.g., \cite{bourbaki}, I.2.7), the differential operators $X^I$ form a basis for the algebra of left-invariant
differential operators in $\G$. Furthermore, let 
$|I|:=i_1+\ldots +i_{n}$ be the order of the differential operator
$X^I$, and let  $d(I):=d_1i_1+\ldots +d_ni_{n}$ be its degree of homogeneity
with respect to group dilations. 
 From the Poincar\'e--Birkhoff-Witt theorem it follows, in particular, that any 
homogeneous linear differential operator in the horizontal  derivatives can be expressed as
a linear combination of the operators $X^I$ of the special form above.

We now recall the notion of convolution in the setting of Carnot groups (see, e.g., \cite{FS}). If $f\in\mc D(\G)$ and
$g\in L^1_{\mathrm{loc}}(\G)$, we set
\begin{equation}\label{group convolution}
f\ast g(x):=\int f(y)g(y^{-1}\cdot x)\,dy\qquad\forall\, x\in \G.
\end{equation}
Furthermore, recall  that  if also  $g$ is a smooth function and $P$
is a left-invariant differential operator, then
$$
P(f\ast g)= f\ast Pg.
$$
  %If $f$ is a real function defined in $\G$, we denote
    %by $\ccheck f$ the function defined by $\ccheck f(p):=
    %f(p^{-1})$, and, if $T\in\mc D'(\G)$, then $\ccheck T$
    %is the distribution defined by $\Scal{\ccheck T}{\phi}
    %:=\Scal{T}{\ccheck\phi}$ for any test function $\phi$.

More generally, we remark that the convolution is   well-defined
whenever $f,g\in\mc D'(\G)$, provided at least one of them
has compact support. In this case,  for any test function $\phi\in\mc D(\G)$, the following identities
hold: 
\begin{equation}\label{convolutions var}
\Scal{f\ast g}{\phi} = \Scal{g}{\ccheck f\ast\phi}
\quad
\mbox{and}
\quad
\Scal{f\ast g}{\phi} = \Scal{f}{\phi\ast \ccheck g}.
\end{equation}

%{  As in \cite{folland_stein},
%we also adopt the following multi-index notation for higher-order derivatives. If 
%$
%I =
%(i_1,\dots,i_{2n+1})
%$ 
%is a multi--index, we set  
%$W^I=W_1^{i_1}\cdots
%W_{2n}^{i_{2n}}\;T^{i_{2n+1}}$. 
%By the Poincar\'e--Birkhoff--Witt theorem, the differential operators $W^I$ form a basis for the algebra of left-invariant
%differential operators in $\G$. 
%Furthermore, if $\G=\he n$, we set 
%$$
%|I|:=i_1+\cdots +i_{2n}+i_{2n+1}
%$$
 %the order of the differential operator
%$W^I$, and   
%$$
%d(I):=i_1+\cdots +i_{2n}+2i_{2n+1}$$
 %its degree of homogeneity
%with respect to group dilations.
%
%

 Suppose now that $f\in\mc E'(\G)$ and $g\in\mc D'(\G)$. If $\psi\in\mathcal D(\G)$, then  it can be shown that
 \begin{equation}\label{convolution by parts}
 \begin{split}
\Scal{(X^If)\ast g}{\psi}&=
 \Scal{X^If}{\psi\ast \ccheck g} =
  (-1)^{|I|}  \Scal{f}{\psi\ast (X^I \,\ccheck g)} \\
&=
 (-1)^{|I|} \Scal{f\ast \ccheck X^I\,\ccheck g}{\psi}.
\end{split}
\end{equation}

\medskip

The following theorem  can be found in  \cite{FS} (see Proposition 1.18).
%(or easily derived from \cite{folland}  \cite{folland_stein}). 

%If $1\le p\le \infty$, $s\ge 0$, we denote by $W^{s,p}(\G)$ the anisotropic Folland-Stein function spaces  (we refer to Section \ref{7.3} below for precise definition and properties).

\begin{theorem}[Hausdorff-Young inequality]  \label{folland cont} 
           If $f\in L^p(\G)$, $g\in L^q(\G)$, $1\le p, q,r \le \infty$, and $\frac1p + \frac1q = 1 + \frac1r$, then $f\ast g\in L^r(\G)$ and
					 $
					\|f\ast g\|_{L^r}\le \|f\|_{L^p}\|g\|_{L^q}\,.
					$ 
         
\end{theorem}

\begin{remark}If $T\in \mc E'(\G)$, and $P$ is a differential operator in $\G$, then $PT\in \mc E'(\G)$, and it turns out that $\supp PT \subset \supp T$
	(see \cite{treves}, Exercise 24.3).
\end{remark}

We collect in the next proposition a few basic properties of the convolution of two  distributions. 

\begin{proposition}\label{treves} The following assertions hold.
\begin{enumerate}
\item If $T\in \mc D'(\G)$ (or, $T\in \mc E'(\G)$, respectively), then the convolution $\phi\mapsto \phi \ast T$ is a continuous linear map of $\mc E(\G)$ (or, $\mc D(\G)$, respectively) into $\mc D(\G)$ (see \cite{treves}, Theorem 27.3).
\item The convolution maps $\mc E(\G)\times \mc D'(\G)$ (or, $\mc D(\G)\times \mc E'(\G)$, respectively) into $\mc D(\G)$ (see \cite{treves}, p. 288).
\item The convolution $(S,T) \mapsto S\ast T$, defined as
$$
\Scal{S\ast T}{\phi}_{\mc D', \mc D} = \Scal{S}{\phi \ast \ccheck T}_{\mc E', \mc E}, 
$$
is a separately continuous bilinear map from $\mc E'(\G)\times \mc D'(\G)$ into $\mc D'(\G)$ (see \cite{treves}, Theorem 27.6).

\end{enumerate}

\end{proposition}

Let $J:\G\longrightarrow\R$ be a mollifier (for the group structure), i.e., $J\in C_c^\infty(\G)$, $J\geq 0$, ${\rm supp}(J)\Subset B(e,1)$, and $\int_{\G}J(x)\,dx=1$.
Note that, if one starts  from a standard mollifier $J$ defined in $(\R, +)$, then the function $J(\|x\|)$ turns out to be a mollifier in $\G$.  Now, given a mollifier $J$, we  define a family of approximations to the identity $\{J_\varepsilon\}_{\varepsilon>0}$ by setting 
$$
J_\varepsilon(x):=\frac{1}{\varepsilon^Q}J(\delta_{1/\varepsilon}x)\,.
$$

We remark explicitly that  $
J_\varepsilon(x)=\ccheck J_\varepsilon(x) 
$  for every $x\in\G$.

Let  $1\le p<+\infty$. If $f\in L^p(\G)$, then $J_\varepsilon\ast f\longrightarrow f$ in $L^p(\G)$ as $\varepsilon\to 0$.
%\st{, and if we denote by $W^{1,p}(\G)$ the  Folland-Stein space, then $J_\varepsilon\ast f\longrightarrow f$ in $W^{1,p}(\G)$ as  $\varepsilon\to 0$}.  
Furthermore, since
$f\ast J_\varepsilon =\ccheck \left(\ccheck J_\varepsilon\ast \ccheck f\right)=\ccheck \left(J_\varepsilon\ast \ccheck f\right)$, the same assertions hold true for $f\ast J_\varepsilon$.

\subsection{Vector Measures in $H\G$ and Riesz Theorem}

Throughout we shall denote by  $C_c(\G, H\G)$ the class of continuous horizontal vector fields  with  compact support in $\G$, and by $C_0(\G, H\G)$ its completion with respect to the uniform norm  
$$\|\Phi\|_\infty=\sup\{|\Phi(x)|_x\,:\, x\in \G\},$$where $\Phi:\G\longrightarrow H\G$.  It turns out that $C_0(\G, H\G)$,  endowed with the uniform norm $\|\cdot\|_\infty$, is a Banach space. Furthermore, since the uniform limit of continuous functions is a continuous function, it follows that $\Phi\in C_0(\G, H\G)$  if, and only if,  $\Phi$ is continuous and for every $\varepsilon>0$ there exists a compact set $\mc K \subset \G$ such that $|\Phi(x)|_x\le \varepsilon$ whenever $x\in \G\setminus \mc K$.

We shall refer to the space $C_0(\G, H\G)$ as the space of  continuous horizontal vector fields  {\it vanishing at infinity}. Exactly as in the Euclidean case,  the linear subspace $\mc D(\G, H\G)$ is dense in $C_0(\G, H\G)$.

%{\color{red} forse si puo' non mettere}
Now we need a substitute for the notion of vector-valued measure in Carnot groups (compare with \cite{MagnaniComi}, Definition 3.5).

Let $\gamma\in\mc M(\G)$ be a Radon measure on $\G$ and let $\alpha:\G\to H\G$ be a (locally) bounded $\gamma$-measurable horizontal vector field.  Hence, there is a naturally defined    linear functional on  $C_c(\G, H\G)$ given by $T_{\alpha\gamma}(\Phi):=\int_\G \scal{\Phi}{\alpha}d\gamma$  (clearly, $T_{\alpha\gamma}$ is
bounded in $C_c(\G, H\G)$ with respect to the $L^\infty$-topology).  As a consequence,  we  can define a notion of {\it vector  measure  $\alpha\gamma$ in $H\G$} by setting
$$C_c(\G, H\G)\ni \Phi\longmapsto 
\int_\G \scal{\Phi}{d(\alpha\gamma)}:=T_{\alpha\gamma}(\Phi).
$$ By density,  this functional extends to a continuous linear functional in $C_0(\G, H\G)$. In the sequel, we shall denote by $\mc M(\G, H\G)$ the space of all vector measures on $\G$ (in the previous sense). 
As previously pointed out, we can write $\alpha=\sum_{i=1}^m\alpha_i X_i$, where the components $\alpha_i:\G\to\R \; (i=1,...,m)$  with respect to the horizontal frame  are now (locally) bounded $\gamma$-measurable functions.
Hence, the vector  measure $\mu=\alpha\gamma$ can be written (in components) as
 $
\mu= (\mu_1,\ldots,\mu_m)=(\alpha_1,\ldots,\alpha_m)\gamma,
$ 
and we get  $$T_{\mu}(\Phi)=\int_\G  \scal{\Phi}{d\mu}=\sum_{i=1}^m\int_\G\Phi_i(x)d\mu_i(x).$$

Since in Carnot groups the horizontal bundle has a global trivialization, we can always argue componentwise. Then it   is not difficult to show that any $T\in  C_0(\G, H\G)^\ast $ can be represented by a  vector measure $\mu$ in $H\G$  as  
$$
T(\Phi)=\int_\G \scal{\Phi}{d\mu}\qquad \forall\,\Phi\in C_0(\G, H\G).
$$  Moreover, due to the density of $\mc D(\G, H\G)$  in $C_0(\G, H\G)$, 
if we take $T\in \mc D(\G, H\G)^\ast$ such that
 $
\sup \left\{T(\Phi)\, :\  \Phi\in \mc D(\G, H\G),\ \|\Phi\|_\infty\le 1\right\}<+\infty,
$ 
we can extend uniquely $T$ to an element of $C_0(\G,H\G)^\ast$. Hence, any $T$ turns out to be associated with a vector measure $\mu \in\mc M(\G, H\G)$. We henceforth set 
$$
\|\mu\|_{\mc M}:=\sup \left\{T(\Phi)\, :\  \Phi\in \mc D(\G, H\G),\ \|\Phi\|_\infty\le 1\right\} =\|T\|_{C_0^\ast}
$$(the symbol $\mc M$  will be  omitted when clear by the context). The identification between the space  $\mc M(\G,H\G)$  of vector measures with finite mass and $ C_0(\G, H\G)^\ast\,$ can be proved using the map $\rho:  \mc M(\G,H\G)\longrightarrow C_0(\G, H\G)^\ast$
	defined by $$\rho(\mu)(\Phi):=\int_\G\langle\Phi, d\mu\rangle=T_\mu(\Phi)\qquad\forall\,\Phi\in C_0(\G, H\G).$$

\subsection{Functions of bounded $H$-variation in Carnot groups}
In this subsection we recall some  known definitions and results concerning functions of  ``intrinsic bounded variation''.

Let $\Omega\subseteq\G$ be an open set. Recall  that a function $f:\Omega \longrightarrow\R$ is said to have  {\it intrinsic  bounded variation in $\Omega$}, and in this case we write $f\in BV_H(\Omega)$, if $f\in L^1(\Omega)$ and
$$
\|D_Hf\|(\Omega):=\sup\left\{\int_\Omega f\,{\rm div}_H \Phi \  dx\,\, :\  \Phi\in \mc D(\Omega, H\Omega),\ \|\Phi\|_\infty\le 1\right\}< +\infty,
$$
where $\|\Phi\|_\infty=\sup\{|\Phi(x)|_x\,:\, x\in \Omega\}$. 

The quantity $
 \|D_Hf\|(\Omega)$ represents the {\it total horizontal variation} (or, {\it $H$-variation}) of the distributional horizontal gradient $D_Hf$ in $\Omega$.

Unless otherwise stated, throughout the paper we shall assume that $\Omega=\G$. In this case, the total $H$-variation of $D_Hf$  in $\G$ will be simply denoted as $\|D_H f\|. $

Note that the  preceding definition can easily be localized. To this aim, let $f\in L^1_{loc}(\Omega)$ and assume that $\|D_H f\|(V)<+\infty$ for every open subset $V\Subset\Omega$.  In this case, we set $f\in BV_{H, loc}(\Omega)$ to denote the space of functions of locally bounded $H$-variation in $\Omega$.

Of course, if $\G$ is commutative and equipped with the Euclidean metric, the previous definitions   coincide  with the classical ones.
There is a wide literature on  $BV_H$-functions in Carnot groups for which we refer, for instance, to  \cite{FSSChouston}, \cite{GN}, \cite{vittone}, and references therein.

By adapting the classical Riesz representation theorem to our setting, one can prove the following ``structure theorem''.

\begin{theorem} If $f\in BV_{H, loc}(\Omega)$, then  $\|D_Hf\|$  is a  Radon measure on $\Omega$. 	 In addition, there exists a bounded $\|D_H f\|$-measurable horizontal section $\sigma_f:\Omega\to H\Omega$ such  that $|\sigma_f(x)|_x =1$ for $\|D_Hf\|$-a.e. $x\in\Omega$, and the following holds
	\begin{equation}\label{riesz}
	\int_\Omega f\,\mathrm{div}_H \Phi \  dx=-\int_\Omega \scal{\Phi}{\sigma_f}\,d\|D_Hf\|\qquad\forall\, \Phi\in \mc D(\Omega, H\Omega).
	\end{equation} 
\end{theorem}

 %Let $f\in BV_{H, loc}(\G)$. We shall denote by $[D_Hf]_a$ the absolutely continuous part of $
%[D_Hf]$, with respect to the Haar measure $\mathscr L^n$ of $\G$, and by $
%[D_Hf]_s$ its singular part. It can be shown that
%$$[D_Hf]=[D_Hf]_a+[D_Hf]_s=\mathcal L^n\res D_H f+[D_Hf]_s 
%$$for any $f\in L^1_{\mathrm loc}(\G)$, so that $D_H f\in L_{\mathrm loc}^1(\G, H\G)$ is the density, with respect to  $\mathscr L^n$, of the absolutely continuous part of $[D_Hf]$.

%\begin{remark}[product rule: a particular case]\label{pr} Let $f\in BV_{H, loc}(\G)$ and $\phi\in\mc D(\G)$. Then, we claim that  $D_H( \phi f)=\phi D_H f  +  f D_H\phi$ in the distributional sense. To prove this claim, we argue componentwise  and take a left-invariant vector field $X_i$ of the horizontal frame $\{X_1,..., X_m\}$. If $\psi\in \mc E(\G)$, then  
%	\begin{eqnarray*} \int_{\G}\left(\phi\,X_i(f)+f X_i\phi\right)\,\psi\,dx=\int_{\G} X_i(f)\,\phi\,\psi\,dx+\int_{\G} X_i(\phi)\,f\,\psi\,dx\\=-\int_{\G}\phi\,f\,X_i(\psi)\,dx=	-\int_{\G}f\left(X_i(\phi\,\psi)-\psi\,X_i\phi\right)\,dx=\int_{\G} X_i(\phi\,f)\,\psi \,dx. 
%	\end{eqnarray*}
%In other words, we have shown that $$X_i(\phi\,f)=\phi\,X_i(f)+f\,X_i\phi\qquad (i=1,...,m)$$ in the distributional sense.  From this, the initial claim easily follows. 
	
%\end{remark}
 
Let $C_H^1(\Omega)$ denote the linear space of  functions $f:\Omega \longrightarrow\R$ such
that the pointwise horizontal partial derivatives $X_1 f, \ldots, X_m f$ are continuous in $\Omega$.

\begin{remark}\label{c1funct} As in the Euclidean case, every function $f\in C_H^1(\Omega)$  belongs to $BV_{H,loc}(\Omega)$.  This follows by integrating by parts. Indeed,  we have $$\int_\Omega f\,{\rm div}_H \Phi \  dx=-\int_{\Omega} \scal{\Phi}{D_H f}\,dx,$$which implies that $\|D_Hf\|(\Omega)=\mathscr L^n\res |D_H f|$, and  $$\sigma_f=\begin{cases}
 \frac{D_H f}{|D_Hf|}\quad \mbox{if}\quad D_H f\neq 0\\ \; \; \; 0\qquad  \,  \mbox{if}\quad D_H f = 0
	\end{cases}\qquad \mathscr L^n\mbox{-\it a.e.}$$
	 
\end{remark}

Let $\Omega=\G$. According to the  previous section's definition, $\mu=\sigma_f \,\|D_Hf\|$ is a vector measure in $H\G$.  Writing $\sigma_f$ with respect to the horizontal frame as 
$
\sigma_f=\sum_{i=1}^m\sigma_{f,i}X_i
$, 
where the components $\sigma_{f,i}: \G\longrightarrow\R\; (i=1,\ldots,m)$ are bounded measurable functions, we have
$\mu=
(\sigma_{f,1},\ldots,\sigma_{f,m})\|D_Hf\|.
$ 
We shall set  $[D_Hf]:=\mu$, so that \eqref{riesz} becomes
\begin{equation}\label{riesz2}
\int_\G f\,\mathrm{div}_H \Phi \  dx=-\int_\G \scal{\Phi}{d[D_Hf]} \qquad\forall\, \Phi\in \mc D(\G, H\Omega).
\end{equation}

%\subsection{Basic properties of $BV_H$-functions}

The following  results are relevant in the theory of bounded $H$-variation functions in Carnot groups (for a proof we refer the reader to the literature quoted above).

 The first one asserts that the (total) $H$-variation is lower semicontinuous with respect to the $L_{loc}^1$-convergence and  follows because the map $f\mapsto \|D_Hf\|(\cdot)$ is the supremum of a family of $L^1$-continuous functionals.
\begin{theorem}\label{lsc} Let $\Omega\subseteq\G$ be an open set. 
Let $\{f_k\}_{k\in\N}$ be a sequence in $BV_H(\Omega)$ such that $f_k\longrightarrow f$ in $L^1_{\mathrm loc}(\Omega)$ as $k\to+\infty$. Then 
$$
\|D_Hf\|(\Omega)\le\liminf_{k\to+\infty}\|D_Hf_k\|(\Omega).
$$
\end{theorem}
The next theorem, in the Euclidean setting, is better known as  the ``Anzellotti-Giaquinta approximation theorem''.
\begin{theorem}\label{AG} Let $\Omega\subseteq\G$ be an open set and  let $f\in BV_H(\Omega)$. Then, there exists a sequence $\{f_k\}_{k\in\N}\subset BV_H(\Omega)\cap C^\infty(\Omega)$ such that $f_k\longrightarrow f$ in $L^1(\Omega)$ as $k\to+\infty$, and
$$
\lim_{k\to+\infty}\|D_Hf_k\|(\Omega)=\|D_Hf\|(\Omega).
$$
\end{theorem}

   If $E\subseteq\G$ is a Borel set, we   set $P_H (E):= \|D_H \chi_E\|$, where $\chi_E$ is the characteristic function of $E$. More generally, if $\Omega\subseteq\G$ is an open set,  we   set $P_H (E, \Omega):= \|D_H \chi_E\|(\Omega)$.  The quantities just defined are the $H$-perimeter of $E$ in $\G$ and in $\Omega$, respectively.

   The next result is the coarea formula for functions of bounded $H$-variation (see, e.g.,   \cite{FSSChouston}, \cite{GN}).

\begin{theorem}[Coarea formula] Let $f\in BV_H(\Omega)$ and set  $E_t:=\{x \in \Omega\,:\; f(x) > t\}$. Then, $E_t$ has finite $H$-perimeter in $\Omega$ for a.e. $t\in\R$ and the following formula  
	holds
	\begin{equation}\label{coarea}
	\|D_H f\|(\Omega) =\int_\R
	P_H (E_t,\Omega) \,dt.
	\end{equation} Conversely, if $f\in L^1(\Omega)$ and $\int_\R
	P_H (E_t,\Omega) \,dt<+\infty$, then $f\in BV_H(\Omega)$.
 
\end{theorem}

%\begin{remark}\label{max rem}
%\st{Let $f\in BV_H(\Omega)$, $t\in\R$, and consider the function $g_t:=\max\{f, t\}$. As in the Euclidean case (see, e.g., \cite{giaquinta_modica_souceck}, p. 340), a useful consequence of the coarea formula is that $g_t\in BV_H(\Omega)$ and that  $\|D_H g_t\|(\Omega)=\|D_H f\|(E_t)\leq \|D_Hf\|(\Omega).$ }
%\end{remark}

Finally, we have to recall a fundamental inequality, whose validity will be of central importance for our next results.

\begin{remark}[Gagliardo-Nirenberg inequality]\label{GNBV}
As is well-known, the classical Gagliardo-Nirenberg inequality has been generalized to  Carnot groups  by many authors (and with different aims); see, e.g., \cite{CDG}, \cite{FGaW}, \cite{FLW_grenoble}, 
\cite{GN}, \cite{GROMOV},  \cite{pansu_thesis}. More precisely, if $f\in \mc D(\G)$,  the inequality states that there exists a ``geometric'' constant $  {\mathscr C}_{_{\! GN}}=\mathscr C_{_{\! GN}}(Q, \G)$ such that
\begin{equation}\label{gagliardo}
\|f\|_{L^{Q/Q-1}}\le \mathscr C_{_{\! GN}}\|D_H f\|_{L^1}.
\end{equation}
The inequality \eqref{gagliardo}  extends to functions in $BV_{H}(\G)$ having compact support. In fact, arguing as in \cite{giusti} (see Theorem 1.28), it is sufficient to approximate $f\in BV_{H}(\G)$ with a sequence $\{f_j\}_{j\in\N}\subset \mc D(\G)$ such that $f_j\longrightarrow f$ in $L^1(\G)$ and $\|D_H f_j\|\longrightarrow \|D_H f\|$ as $j\to+\infty$.  
Then, by  \eqref{gagliardo} the sequence is uniformly bounded in the $L^{Q/Q-1}$-norm and hence there exists a subsequence weakly convergent to some $f_0\in L^{Q/Q-1}(\G)$. But since  $f_j\longrightarrow f$ in $L^1(\G)$ as $j\to+\infty$, it follows  that  $f_j\rightharpoonup f=f_0$ in $L^{Q/Q-1}(\G)$ as $j\to+\infty$ and the proof is achieved by using the weak lower semicontinuity of the $L^{Q/Q-1}$-norm (see, e.g., \cite{Brezis}, Proposition 3.5).

\end{remark}

\section{The space $BV^{{{Q}/{Q-1}}}_H(\G)$}
We  introduce another  intrinsic $BV_H$-type  space, which is in fact a subspace of $ L^{{{Q}/{Q-1}}}(\G)$, where $Q$ denotes the homogeneous dimension (equal to the Hausdorff dimension) of  $\G$; see \eqref{dim_omo}. In the Euclidean setting this space was introduced and studied by De Pauw and Torres in \cite{DePauw-Torres}. 
\begin{definition}
The space $BV^{{{Q}/{Q-1}}}_H(\G)$ is the set of functions  $f\in L^{{{Q}/{Q-1}}}(\G)$ whose distributional gradient $D_Hf$ is a  finite vector measure, i.e., 
$$
\|D_Hf\|:=\|D_Hf\|(\G)=\sup\left\{\int_\G f\,{\rm div}_H \Phi \  dx\,\, :\  \Phi\in \mc D(\G, H\G),\ \|\Phi\|_\infty\le 1\right\}< +\infty.
$$
\end{definition}
The space $BV^{{{Q}/{Q-1}}}_H(\G)$ is a Banach space when endowed with the norm 
$$\|f\|_{L^{{{Q}/{Q-1}}}}+\|D_Hf\|.
$$

 Note also that $BV^{{{Q}/{Q-1}}}_H(\G)\subset BV_{H,{loc}}(\G)$.

The next result shows the  lower semicontinuity  of the $H$-variation with respect to the weak convergence in $L^{{{Q}/{Q-1}}}(\G)$.
\begin{theorem}\label{lsc seconda}
Let $\{f_k\}_{k\in\N}$ be a sequence in $BV^{{{Q}/{Q-1}}}_H(\G)$ such that $f_k\rightharpoonup f$ in $L^{{{Q}/{Q-1}}}(\G)$ as $k\to+\infty$. Then
$$
\|D_Hf\|\le\liminf_{k\to+\infty}\|D_Hf_k\|.
$$
\end{theorem}
\begin{proof}
We consider the functional $\int_\G f\,{\rm div}_H \Phi \  dx\,$ with $\Phi\in \mc D(\G, H\G)$ and $\|\Phi\|_\infty\le 1$.
Since ${\rm div}_H \Phi \in L^Q(\G)$ and $f_k\rightharpoonup f$ in $L^{{{Q}/{Q-1}}}(\G)$ as $k\to+\infty$,  we have
$$
\int_\G f\,{\rm div}_H \Phi \  dx\,=\lim_{k\to +\infty}\int_\G f_k\,{\rm div}_H \Phi \  dx\,.
$$
By assumption,  $\{f_k\}_{k\in\N}\subset BV^{{{Q}/{Q-1}}}_H(\G)$, and hence $\int_\G f_k\,{\rm div}_H \Phi \  dx\le \|D_Hf_k\|$. Thus
$$
\int_\G f\,{\rm div}_H \Phi\  dx\,\le\liminf_{k\to+\infty}\|D_Hf_k\|,
$$
and the conclusion follows by taking the supremum on the left-hand side over all $\Phi$ in $\mc D(\G, H\G)$ such that $\|\Phi\|_\infty\le 1$.
\end{proof}

\subsection{An approximation result for $BV^{{{Q}/{Q-1}}}_H(\G)$}

%\medskip
%
%The vector measure $D_Hf$ can be convolved with a Friederichs mollifier
%
%\medskip

We start with an approximation result that yields as corollaries a Gagliardo-Nirenberg inequality for functions in $BV^{{{Q}/{Q-1}}}_H(\G)$ and  a compactness result in $BV^{{{Q}/{Q-1}}}_H(\G)$. The results in this subsection generalize the corresponding Euclidean ones in \cite{DePauw-Torres}.
\begin{theorem}\label{density}
Let $f\in BV^{{{Q}/{Q-1}}}_H(\G)$. Then, there exists a sequence $\{f_j\}_{j\in \N}\subset\mc D(\G)$ such that: 
\begin{itemize}
	\item[(i)] $f_j\rightharpoonup f$ in $L^{{{Q}/{Q-1}}}(\G)$ as $j\to+\infty$ and  $\sup_j \|D_Hf_j\|<+\infty$.\\

In addition, the sequence $\{f_j\}_{j\in \N}$ satisfies: \\
	\item[(ii)]$\lim_{j\to+\infty}\|D_Hf_j\|_{L^1}=\|D_Hf\|$.
\end{itemize} 
\end{theorem}

\begin{proof}

The proof is divided in several steps.
\begin{itemize}
	\item[\bf Step 1.] Consider a family of approximations to the identity $\{J_\varepsilon\}_{\varepsilon>0}$  (see Section 2) and remember that $J_\varepsilon=\ccheck J_\varepsilon$.  Since $J_\varepsilon\ast f\longrightarrow f$ in $L^{{{Q}/{Q-1}}}(\G)$ as $\varepsilon\to 0^+$, one has obviously $J_\varepsilon\ast f\rightharpoonup f$ in $L^{{{Q}/{Q-1}}}(\G)$ as  $\varepsilon\to 0^+$.
	In addition, it follows from \eqref{convolutions var} that if $\Phi\in \mc D(\G, H\G)$ and $\|\Phi\|_\infty\le 1$, then
	$$
	\Scal{J_\varepsilon\ast f}{X_i\Phi} = \Scal{f}{\ccheck J_\varepsilon\ast X_i\Phi}=\Scal{f}{ J_\varepsilon\ast X_i\Phi}=\Scal{f}{ X_i(J_\varepsilon\ast\Phi)}.
	$$
Hence 
	$$
	\int_\G (J_\varepsilon\ast f)\,\divg \Phi\,dx=\int_\G  f\,\divg(J_\varepsilon\ast\Phi)\,dx.
	$$
Now since $\|J_\varepsilon\ast\Phi\|_\infty\le \|\Phi\|_\infty\le 1$, taking the supremum on the right-hand side,  we get
	$$
	\int_\G (J_\varepsilon\ast f)\,\divg \Phi\, dx\le \|D_H f\|.
	$$
In turn, since $J_\varepsilon\ast f\in C^\infty(\G)$, taking  the supremum on the left-hand side over all $\Phi\in \mc D(\G, H\G)$ such that $\|\Phi\|_\infty\le 1$, we obtain
\begin{equation}\label{stima hvar}
	\|D_H (J_\varepsilon\ast f)\|_{L^1}=	\|D_H (J_\varepsilon\ast f)\| \le \|D_H f\| 
\end{equation}
	for every $\varepsilon>0$; see, e.g., Remark \ref{c1funct}. So let $\{\varepsilon_k\}_{k\in\N}$ be a strictly decreasing sequence such that $\varepsilon_k\to 0$ as $k\to+\infty$.
Using the lower semicontinuity property  in Theorem \ref{lsc seconda} together with \eqref{stima hvar}, it follows eventually that
\begin{equation}\label{22}	\lim_{k\to +\infty}\|D_H (J_{\varepsilon_k}\ast f)\|_{L^1} =\|D_H f\|. 
\end{equation}

	\item[\bf Step 2.] Starting from \eqref{22}, it is clear that there must exist a subsequence $\big\{J_{\varepsilon_{{k}_j}}\ast f\big\}_{j\in\N}$ of $\{J_{\varepsilon_{k}}\ast f\}_{k\in\N}$ such that \begin{equation}
	\label{22bis}\|D_H (J_{\varepsilon_{k_j}}\ast f)\|_{L^1} \le \|D_H f\|+ \frac{1}{j}\qquad \forall\,j\in\N.
	\end{equation}

		\item[\bf Step 3.] Let us fix a sequence of cut-off functions $\{g_i\}_{i\in\N}\subset\mc D(\G)$  such that for any $i\in\N$   ${\rm supp} (g_i)\subset {B(e,2i)}$, $g_i\equiv 1$ in  ${B(e,i)}$, and $ \sup_i\|D_H g_i\|<+\infty.$
We   have 
	\begin{equation}\label{product}
	D_H((J_{\varepsilon_{k_j}}\ast f)g_i)=g_i D_H (J_{\varepsilon_{k_j}}\ast f)+ (J_{\varepsilon_{k_j}}\ast f) D_H g_i.
	\end{equation}
	%\st{and it is immediate to see that $(J_{\varepsilon_{k_j}}\ast f) g_i\in BV_H(\G)\subset BV^{{{Q}/{Q-1}}}_H(\G)$.} 
	
Let us start by estimating the second term of the right hand side above. Let $j\in\N$ be fixed. Since $J_{\varepsilon_{k_j}}\ast f\in L^{{{Q}/{Q-1}}}(\G)$, it follows  that
	
\begin{eqnarray*} \limsup_{i\to+\infty}\int_{\G}&\!\!&\left| (J_{\varepsilon_{k_j}}\ast f ) D_H g_i\right|\,dx=\limsup_{i\to+\infty}\int_{\G\setminus B(e,i)}\left| (J_{\varepsilon_{k_j}}\ast f ) D_H g_i\right|\,dx\\&\le & \limsup_{i\to+\infty}\left(\int_{\G\setminus B(e,i)}\left| (J_{\varepsilon_{k_j}}\ast f ) \right|^{{{Q}/{Q-1}}}\,dx\right)^{{{Q-1}/{Q}}} \|D_H g_i\|_{L^Q}\\&=& 0.
\end{eqnarray*}
	With this estimate in mind, and by means of \eqref{22bis},  it can be shown that there exists a strictly increasing sequence $\{i_j\}_{j\in\N}$ such that   \begin{eqnarray}\nonumber \int_{\G}\left| D_H \left((J_{\varepsilon_{k_j}}\ast f ) g_i\right)\right|\,dx&\le& \int_{\G}\left| D_H  (J_{\varepsilon_{k_j}}\ast f)  \right|\,dx+ \frac{1}{j}\\\label{eq inter} \\\nonumber &\le& \|D_H f\|+ \frac{2}{j}\qquad\forall\,j\in\N.\end{eqnarray}

\item[\bf Step 4.]  
Let us set $$f_j:= (J_{\varepsilon_{k_j}}\ast f) g_{i_j}\qquad \forall\,j\in\N.$$ From Step 3 it follows in particular  that $\sup_j \|D_Hf_j\|<+\infty$.
Let us to show that  $f_j\rightharpoonup f$ in $L^{{{Q}/{Q-1}}}(\G)$ as $j\to+\infty$. If we take $g\in L^Q(\G)$, we have 

 \begin{eqnarray*} \left|\int_{\G}g\left(f- (J_{\varepsilon_{k_j}}\ast f)  g_{i_j} \right)\,dx\right| \leq \int_{\G}|g|\left|f- (J_{\varepsilon_{k_j}}\ast f)  \right|\,dx + \int_{\G}|g|\left| J_{\varepsilon_{k_j}}\ast f \right||1-g_{i_j}|\,dx\\\leq  \|g\|_{L^Q} \| f- (J_{\varepsilon_{k_j}}\ast f) \|_{L^{{{Q}/{Q-1}}}}+ \left(\int_{\G\setminus B(e, 2i_j)} |g|^Q \,dx\right)^{1/Q} \|f\|_{L^{{{Q}/{Q-1}}}} .\end{eqnarray*}
Since both addends of the right-hand side vanish as  $j\to+\infty$, assertion (i) is proved.  
Finally, using the inequalities \eqref{eq inter} together with the lower semicontinuity property in Theorem \ref{lsc seconda}, it follows that $\lim_{j\to+\infty}\|D_Hf_j\|_{L^1}=\|D_Hf\|$, which proves  (ii).

\end{itemize}

\end{proof}

%\begin{remark}\label{rem appr} Let  $f\in BV^{{{Q}/{Q-1}}}_H(\G)$ and $\phi\in\mc D(\G)$. Arguing as in Step 1 of the preceding proof, we get in particular that  $f\ast \phi \in  BV^{{{Q}/{Q-1}}}_H(\G)$ and that \begin{equation}\label{stima hvar-bis}
	%\|D_H (f\ast\phi)\|_{L^1}\le \|D_H f\|\|\phi\|_{L^1}.
	%\end{equation}	\end{remark}

\begin{corollary}[Gagliardo-Nirenberg inequality in $BV_H^{Q/Q-1}(\G)$] \label{gn bv-type}  Let $f\in BV_H^{Q/Q-1}(\G)$. Then \begin{equation} \label{gagliardo-bis}
	\|f\|_{L^{Q/Q-1}}\leq {\mathscr C}_{_{\! GN}}\|D_H f\|.\end{equation}
\end{corollary}

\begin{proof} The proof follows by approximating $f$ as in  Theorem \ref{density}, using  inequality \eqref{gagliardo} for functions in $\mc D(\G)$, and then applying the
 weak lower semicontinuity of the $L^{Q/Q-1}$-norm.

\end{proof}

  \begin{remark}\label{oss gn} Let $f\in BV^{{{Q}/{Q-1}}}_H(\G)$. By \eqref{gagliardo-bis} it follows that the $H$-variation
  	$\|D_H f\|$   is an equivalent norm to $\|f\|_{L^{{{Q}/{Q-1}}}}+\|D_Hf\|$. For this reason, in the sequel  the $H$-variation  will be taken as a norm and we shall set $$\|f\|_{BV_H^{{Q}/{Q-1}}}:=\|D_Hf\|.$$ 
  Note also that \eqref{gagliardo-bis} immediately implies the continuous embedding
  	\begin{equation}\label{immersione}
  	BV_H(\G)  \hookrightarrow BV^{{{Q}/{Q-1}}}_H(\G).
  	\end{equation}
  \end{remark}

As a corollary of Theorem \ref{density} and of the Gagliardo-Nirenberg inequality, we obtain the following
 compactness result.
\begin{corollary}[compactness]\label{compact}
  Let $\{f_k\}_{k\in\N}$ be a  sequence in $BV_H^{Q/Q-1}(\G)$ satisfying
	$$\sup_k\|D_Hf_k\|<+\infty. $$ Then, there exists a subsequence $\{f_{k_j}\}_{j\in\N}$ and a function $f\in BV_H^{Q/Q-1}(\G)$ such that
	$$f_{k_j}\rightharpoonup  f\quad  \mathit{in}\quad L^{Q/Q-1}(\G) \quad \mathit{as}\quad j\to+\infty.$$
\end{corollary}. 
\begin{proof}
Since $\sup_k\|D_Hf_k\|<+\infty $, by Corollary \ref{gn bv-type} $\{f_k\}_{k\in\N}$ is  equibounded in $L^{Q/Q-1}(\G)$. Hence 
there exists a subsequence  $\{f_{k_j}\}_{j\in\N}$ that weakly converges {in}  $L^{Q/Q-1}(\G)$ to some function $f$ (see, e.g. \cite{Brezis}, Theorem 3.18). By Theorem \ref{lsc seconda}, $\|D_Hf\|\le \liminf_{j\to +\infty}\|D_Hf_{k_j}\|$. Thus, using the equiboundeness of $\|D_Hf_{k_j}\|$, it follows that $f\in BV_H^{Q/Q-1}(\G)$.
\end{proof}

\section{Charges vanishing at infinity}
In this section we shall define a subspace of $\big(BV^{{{Q}/{Q-1}}}_H(\G)\big)^\ast$, denoted by ${\bf Ch}_0(\G)$, and we shall  investigate the relationship between its dual and the space $\big(BV^{{{Q}/{Q-1}}}_H(\G)\big)$. 

The results of this section will be used  later, in order to  define  a divergence-type operator from $ C_0(\G, H\G)$ to $ {\bf Ch}_0(\G)$, which will turn out to be a bounded linear operator.

 In rough terms, this  operator will be the right substitute for the horizontal divergence operator $\mathrm{div}_H$, when acting on $C_0(\G, H\G)$, and we shall prove that is a surjective operator, which means that we can find a solution in $ C_0(\G, H\G)$ to the equation  $\mathrm{div}_H\Phi=F$, whenever $F\in {\bf Ch}_0(\G)$.

The presentation and results in this section are largely inspired by those in \cite{DePauw-Torres}. 
\begin{definition}\label{def BVstar}
Given a sequence $\{f_j\}_{j\in\N}$ in $ BV^{{{Q}/{Q-1}}}_H(\G)$ we write $$f_j \twoheadrightarrow 0\qquad (j\to+\infty)$$
if and only if $f_j \rightharpoonup 0$ in $L^{Q/Q-1}(\G)$ as $j\to+\infty$ and $\sup_j \|D_H f_j\| <+\infty$.
%\begin{itemize}
%\item[(i)] ;
%\item[(ii)] 
%\end{itemize}
\end{definition}
 More generally, if $f\in BV^{{{Q}/{Q-1}}}_H(\G)$, we  write $f_j - f\twoheadrightarrow 0$  as $j\to+\infty$  whenever $f_j \rightharpoonup f$ in $L^{Q/Q-1}(\G)$ as $j\to+\infty$ and $\sup_j \|D_H f_j\| <+\infty$. 
 	
%{\color{red} LEVEREI TUTTO QUESTO	 In fact, note that  by the triangle inequality we have $$\|D_H(f-f_j)\|\leq \|D_H f_j\|+\|D_Hf\|<+\infty\qquad\forall\,j\in\N.$$
%As a consequence, we can reformulate Theorem \ref{density} by saying that for every $f\in BV^{{{Q}/{Q-1}}}_H(\G)$ there exists $\{f_j\}_{j\in\N}\subset\mc D(\G)$ such that $f_j \twoheadrightarrow f$ and $\|D_H f_j\|\longrightarrow\|D_H f\|$ as $j\to+\infty$. }
  	
% \end{remark}

\begin{definition}[Charges vanishing at $\infty$] Let  $F: BV^{{{Q}/{Q-1}}}_H(\G)\longrightarrow \R $  be a linear functional. We say that $F$ is a {\rm charge vanishing at $\infty$} if and only if  $$\Scal{ F }{f_j}\xrightarrow[j\longrightarrow+\infty]{} 0$$ for any sequence $\{f_j\}_{j\in\N}\subset BV^{{{Q}/{Q-1}}}_H(\G)$ such that $f_j \twoheadrightarrow 0$ as $j\to+\infty$.

\end{definition}

From now on we shall denote by ${\bf Ch}_0(\G)$ the class of all charges vanishing at $\infty$.\\
\begin{remark}
It is clear that ${\bf Ch}_0(\G)$ is a (real)  vector space. We set $$\|F\|_{{\bf Ch}_0}:=\sup\left\{ \Scal{ F }{f}\;:\; f\in BV^{{{Q}/{Q-1}}}_H(\G),\;\; \|D_H f\|\leq 1 \right\}.$$
 Notice that $\|F\|_{{\bf Ch}_0}<+\infty$ whenever $F\in {\bf Ch}_0(\G)$.   In fact, there exists a sequence $\{f_j\}_{j\in\N}\subset BV^{{{Q}/{Q-1}}}_H(\G)$ with $\|D_H f_j\|\leq 1$ such that $\Scal{ F }{f_j}\longrightarrow \|F\|_{{\bf Ch}_0}$ as   $j\to+\infty$. 
	By Proposition \ref{compact}, there exist $f\in BV^{{{Q}/{Q-1}}}_H(\G) $ and a subsequence $\{f_{j_k}\}_{k\in\N}$  such that $f_{j_k}-f\twoheadrightarrow 0$ as   $k\to +\infty$. As a consequence,  $\Scal{ F }{f_{j_k}-f}\longrightarrow 0 $ as $k\to+\infty$.
 Thus $$\Scal{ F }{f}=\lim_{k\to+\infty}\Scal{ F }{f_{j_k}}=\|F\|_{{\bf Ch}_0}<+\infty. $$	
\end{remark}

From this remark it follows that $\|\cdot\|_{{\bf Ch}_0}$ is a norm on ${{\bf Ch}_0}(\G)$.  We also observe that  ${{\bf Ch}_0}(\G)\subset \big(BV^{{{Q}/{Q-1}}}_H(\G)\big)^\ast$ 
and that for any $F\in {\bf Ch}_0(\G)$ we have $$\|F\|_{{\bf Ch}_0}=\|F\|_{\big(BV^{{{Q}/{Q-1}}}_H(\G)\big)^\ast}.$$

\begin{proposition}
	The  space ${{\bf Ch}_0}(\G)$ is a Banach space under the norm  $\|\cdot\|_{{\bf Ch}_0}$.
\end{proposition}

\begin{proof} We show that each Cauchy sequence $\{F_k\}_{k\in\N}\subset {{\bf Ch}_0}(\G)$ converges to an element of ${{\bf Ch}_0}(\G)$.  To this end, note  that $\{F_k\}_{k\in\N}$ has to converge to some $F\in  \big(BV^{{{Q}/{Q-1}}}_H(\G)\big)^\ast$, hence for any $\varepsilon >0$ there exists $k_\varepsilon\in \N$ such that  $\|F-F_k\|_{\big(BV^{{{Q}/{Q-1}}}_H(\G)\big)^\ast}<\varepsilon$ for any $k>k_\varepsilon$.

Let now $k>k_\varepsilon$
	and let $\{f_j\}_{j\in\N}\subset BV^{{{Q}/{Q-1}}}_H(\G)$ be any sequence such that $f_{j}\twoheadrightarrow 0$ as  $j\to +\infty$. Furthermore, set $$ \mathscr K:=\sup_j \|D_H f_j\|.$$ 
	%Clearly, for any $\epsilon>0$ there exists an index $k=k_\epsilon$ such that $$\|F-F_k\|_{\big(BV^{{{Q}/{Q-1}}}_H(\G)\big)^\ast}\leq \epsilon.$$
	For every $j\in\N$ 
	\begin{eqnarray*} |\Scal{F}{f_j}| &\leq&|\Scal{F-F_k}{f_j}|+ |\Scal{F_k}{f_j}|\\&\leq&   \mathscr  K \,\| {F-F_k}\|_{\big(BV^{{{Q}/{Q-1}}}_H(\G)\big)^\ast}+ |\Scal{F_k}{f_j}| \\&\leq&   \mathscr K \,\varepsilon+ |\Scal{F_k}{f_j}|.
 	\end{eqnarray*}In turn, this implies that $$\limsup_{j\to+\infty} |\Scal{F}{f_j}|\leq    \mathscr K \,\varepsilon.$$ From the arbitrariness of $\varepsilon>0$ we get that $F\in {\bf Ch}_0(\G)$. 
	\end{proof}

	\subsection{An example of charge vanishing at $\infty$}
%We now introduce  two important examples of charges  
	
	Since $BV^{{{Q}/{Q-1}}}_H(\G)\subset L^{{{Q}/{Q-1}}}(\G)$, we can state the following definition.
	\begin{definition} For any $f\in L^Q(\G)$, let $\Lambda(f):BV^{{{Q}/{Q-1}}}_H(\G)\longrightarrow\R$  be the linear functional  defined by $$\Scal{\Lambda(f)}{g}:=\int_{\G}f g\,dx.$$
	\end{definition}

%Notice that  $\Lambda(f)$ is defined on  the space $BV^{{{Q}/{Q-1}}}_H(\G)$, which is a subspace of $L^{{{Q}/{Q-1}}}(\G)$.

\begin{proposition}\label{lambda operator} If $f\in L^Q(\G)$, then $\Lambda(f)\in {\bf Ch}_0(\G)$ and  $\|\Lambda(f)\|_{{\bf Ch}_0}\leq {\mathscr C}_{_{\! GN}} \|f\|_{L^Q}$. 
	 Thus, the  linear operator  $\Lambda: L^Q(\G)\longrightarrow {\bf Ch}_0(\G)$  is a bounded linear operator whose  norm is bounded  by the Gagliardo-Nirenberg constant ${\mathscr C}_{_{\! GN}}$. 
\end{proposition}
	\begin{proof}
		Let  $\{g_j\}_{j\in\N}\subset BV^{{{Q}/{Q-1}}}_H(\G)$ be a sequence such that $g_{j}\twoheadrightarrow 0$ as  $j\to +\infty$.\\ In particular, this sequence weakly converges to $0$ in  $L^{{{Q}/{Q-1}}}(\G)$. So we get that $$ {\Scal{\Lambda(f)}{g_j}} = \int_{\G}f g_j\,dx\xrightarrow[j\to+\infty]{ } 0,$$which shows that $\Lambda(f) \in{\bf Ch}_0(\G)$.
	 Moreover, for any $g\in BV^{{{Q}/{Q-1}}}_H(\G)$ we have $$\left|\Scal{\Lambda(f)}{g}\right|\leq \|f\|_{L^Q}\,\|g\|_{L^{{{Q}/{Q-1}}}}\leq {\mathscr C}_{_{\! GN}}\|D_H g\|\,\|f\|_{L^Q},$$ where we have used  H\"older inequality and the Gagliardo-Nirenberg  inequality \eqref{gagliardo-bis}. Hence $$\|\Lambda(f)\|_{{\bf Ch}_0}\leq {\mathscr C}_{_{\! GN}} \|f\|_{L^Q}.$$
	\end{proof}
	
	We would like to show that the image $\mc R(\Lambda)$  of $\Lambda$ is dense in ${\bf Ch}_0(\G)$  or, equivalently, that any charge vanishing at infinity can be approximated by  a charge in $\mc R(\Lambda)$.

%\subsection{Approximation of  charges vanishing at infinity}

As already recalled in Section \ref{preliminar} (see, e.g., Proposition \ref{treves}), we notice that in distribution theory 
	the common way to define the convolution between a distribution $F$ and a test function $\phi$ is as follows: $$\Scal{F\ast\phi}{\psi}:=\Scal{F}{\psi\ast\ccheck{\phi}}_{\mc D', \mc D} \qquad\forall\;\psi\in\mc D(\G).$$ Now, let $F\in {{\bf Ch}_0}(\G)$ and  $\phi\in\mc D(\G)$: our aim is to define a new  charge $F\ast\phi$.

	More precisely, let $g\in BV^{{{Q}/{Q-1}}}_H(\G)$
 and $\phi\in\mc D(\G)$. Arguing as in Step 1 of the  proof of Proposition \ref{density}, we get  that  $g\ast \ccheck\phi \in  BV^{{{Q}/{Q-1}}}_H(\G)$ and that \begin{equation}\label{stima hvar-bis}
	\|D_H (f\ast\phi)\|_{L^1}\le \|D_H f\|\|\phi\|_{L^1}.
	\end{equation}
	%Then, from Remark \ref{rem appr} 
	%we immediately get that $g\ast \ccheck\phi\in BV^{{{Q}/{Q-1}}}_H(\G) $.  
	This  motivates the following definition.
	
	\begin{definition}
	Let  $F\in {{\bf Ch}_0}(\G)$ and  $\phi\in\mc D(\G)$. We define  the linear functional  $$F\ast\phi:BV^{{{Q}/{Q-1}}}_H(\G)\longrightarrow\R$$  by setting $$ BV^{{{Q}/{Q-1}}}_H(\G) \ni g\longmapsto \Scal{F\ast\phi}{g}:=\Scal{F}{g\ast\ccheck\phi}.$$
	\end{definition}
	
	%To prove  the next two propositions we follow almost verbatim the corresponding ones in the paper \cite{DePauw-Torres}, hence we could omit the proofs but we prefer to  rewrite them for sake of completeness.
	\begin{proposition}Let  $F\in {{\bf Ch}_0}(\G)$ and  $\phi\in\mc D(\G)$. Then $F\ast\phi\in{\bf Ch}_0(\G)\cap\mc R(\Lambda)$.
		
	\end{proposition}

\begin{proof} The proof follows almost verbatim the corresponding one in  \cite{DePauw-Torres}, Proposition 4.1, and we  sketch it for the reader's convenience. When one restricts $F$ to $\mc D(\G)$, the restricted functional is a distribution. Thus, the convolution $F\ast\phi$ is a well-defined distribution, which is actually a smooth function. Thus, there must exist $f\in C^\infty(\G)$ such that 
\begin{equation}\label{eq lambda op}
 \Scal{F\ast\phi}{\psi}= \int_{\G} f\psi\,dx\qquad\forall\psi\in\mc D(\G).
\end{equation} 	
	
Moreover, it is not difficult to see that the function $f$ belongs to $L^Q(\G)$. In fact, let $\{\psi_j\}_{j\in\N}\subset\mc D(\G)$ be a sequence such that $\psi_j\rightharpoonup 0$ in $L^{Q/Q-1}$ and  $\|\psi_j\|_{L^{{{Q}/{Q-1}}}}\xrightarrow[j\to+\infty]{}0$. 
 \\We clearly have the following:
\begin{eqnarray*}\sup_j\|D_H(\psi_j\ast\ccheck \phi)\|&=&\sup_j\|D_H(\psi_j\ast\ccheck \phi)\|_{L^1}=\sup_j\|\psi_j\ast D_H(\ccheck \phi)\|_{L^1}\\&\leq&\sup_j\|\psi_j\|_{L^{{{Q}/{Q-1}}}}\|D_H(\ccheck\phi)\|_{L^{{{Q}/{Q+1}}}}\\&=&\sup_j\|\psi_j\|_{L^{{{Q}/{Q-1}}}}\|D_H \phi\|_{L^{{{Q}/{Q+1}}}}<+\infty,\end{eqnarray*}where we have used the Hausdorff-Young inequality (see  Theorem \ref{folland cont}). Moreover, for any $g\in L^Q(\G)$ we have $$\int_{\G}{g(\psi_j\ast\phi)}\,dx=\int_{\G}{\psi_j(g\ast\phi)}\,dx.$$ Since $g\ast\phi\in L^{Q}(\G)$  and  $\psi_j\rightharpoonup 0$ weakly in $L^{{Q}/{Q-1}}(\G)$  as $j\to+\infty$,  the right-hand side of the last equality tends to $0$ as $j\to+\infty$. In particular,  this implies that  $\psi_j\ast\phi\twoheadrightarrow  0$ as $j\to+\infty$ and that $$\Scal{F}{\psi_j\ast\phi}= \Scal{F\ast\phi}{\psi_j}\xrightarrow[j\to+\infty]{}0.$$ Thus, the linear functional $F\ast\phi$ turns out to be continuous in $\mc D(\G)$ (with respect to the  topology of $L^{{Q}/{Q-1}}(\G)$). The density of $\mc D(\G)$ in $L^{{Q}/{Q-1}}(\G)$ implies that $F\ast\phi$ can be uniquely extended to a bounded  linear functional on $L^{{Q}/{Q-1}}(\G)$. Thus, it follows from the Riesz representation theorem that $f\in L^Q(\G)$.  

Note that  since $f\in L^Q(\G)$,  Proposition \ref{lambda operator}  implies that  $\Lambda(f)\in{\bf Ch}_0(\G)$.

 We are left to show  that $\Lambda(f)=F\ast \phi$, which means that $F\ast \phi\in \mc R(\Lambda)$.

In fact, this is equivalent to show that equation \eqref{eq lambda op} holds true whenever $\psi\in BV^{{{Q}/{Q-1}}}_H(\G)$.

\noindent By  {Theorem} \ref{density}
% (see also Remark \ref{rem appr mg})
 we  can take a sequence $\{\psi_j\}_{j\in\N}\subset\mc D(\G)$ such that $\psi_j \twoheadrightarrow \psi$ as $j\to+\infty$. Hence, from \eqref{eq lambda op} we get that $\Scal{F\ast\phi}{\psi_j} = \int_{\G} f\psi_j\,dx$ for every $j\in\N$ and $$\int_{\G} f\psi_j\,dx\xrightarrow[j\to+\infty]{}\int_{\G} f\psi \,dx.$$We also observe that
$$\Scal{F\ast\phi}{\psi_j} =\Scal{F}{\psi_j\ast\ccheck\phi} \xrightarrow[j\to+\infty]{} \Scal{F}{\psi\ast\ccheck\phi} =\Scal{F\ast\phi}{\psi},$$which is true because $F\in {{\bf Ch}_0}(\G)$ and  $\psi_j\ast\ccheck\phi \twoheadrightarrow \psi\ast\ccheck\phi$ as $j\to+\infty$.
As a consequence  $$ \Scal{F\ast\phi}{\psi}= \int_{\G} f\psi \,dx\qquad\forall\,\psi\in  BV^{{{Q}/{Q-1}}}_H(\G),$$as wished.

	 \end{proof}

	% Remark \ref{approxtoid}

Let $\{J_\varepsilon\}_{\varepsilon>0}$ be  a family of approximations to the identity associated with a symmetric kernel (i.e.,  $J_\varepsilon(x)=\ccheck J_\varepsilon(x)$ for every $x\in\G$). Let $\{\varepsilon_k\}_{k\in\N}$ be a strictly decreasing sequence such that $\varepsilon_k\to 0^+$ as $k\to+\infty$.

		\begin{proposition}  
			
			Let  $F\in {{\bf Ch}_0}(\G)$ and let $\{J_{\varepsilon_k}\}_{k\in\N}\subset\mc D(\G)$ be as above. Then $$\|F- F\ast J_{\varepsilon_k}\|_{{\bf Ch}_0}\xrightarrow[k\to+\infty]{}0.$$
			
		\end{proposition}
		 We  omit this proof since  it looks very similar to the corresponding one in \cite{DePauw-Torres} (see Proposition  4.2).

	\begin{remark}\label{density2}
		An immediate consequence of this approximation result is the density of $\mc R(\Lambda)$ in the space ${\bf Ch}_0(\G)$ of all charges vanishing  at $\infty$.
	\end{remark}

	\begin{remark}\label{=bis} For any $\Phi\in \mc D(\G, H\G)$ with $\|\Phi\|_\infty\le 1$, let us consider the charge $\Lambda({\rm div}_H\Phi)$. Since$$\Scal{\Lambda({\rm div}_H\Phi)}{g}= \int_{\G}g \,{\rm div}_H\Phi \,dx\qquad\forall\,g\in  BV^{{{Q}/{Q-1}}}_H(\G),$$ 
	we infer that	$
	\Scal{\Lambda({\rm div}_H\Phi)}{g}\le \|D_Hg\|\,.
	$
	Thus, if $g\in   BV^{{{Q}/{Q-1}}}_H(\G)$ and $\|D_Hg\|\le 1$, we immediately get that
	\begin{equation}\label{norma_lambda}
		\|{\Lambda({\rm div}_H\Phi)}\|_{{\bf Ch}_0}\le 1\,.
	\end{equation}
	
	\end{remark}

	\begin{proposition}\label{isom} There exists  a linear bijective operator ${\mathit ev}: BV^{{{Q}/{Q-1}}}_H(\G)\longrightarrow {\bf Ch}^\ast_0(\G)$, given by  $$\Scal{{\mathit ev}(f)}{F} := \Scal{F}{f}\qquad\forall\,f\in BV^{{{Q}/{Q-1}}}_H(\G)\quad \forall\, F\in {\bf Ch}_0(\G).$$ 
	\end{proposition}
	
	\begin{proof}  It is obvious that $\mathit  ev$ is a linear operator.  Furthermore, since $$|\Scal{{\mathit ev}(f)}{F}|= |\Scal{F}{f}|\leq \|F\|_{{\bf Ch}_0}\|D_Hf\|,$$  it follows that the operator ${\mathit ev}$ maps $BV^{{{Q}/{Q-1}}}_H(\G)$ onto ${\bf Ch}^\ast_0(\G)$. In order to show that $\it ev$  is also injective,  let   $f\in BV^{{{Q}/{Q-1}}}_H(\G)$ be such that ${\mathit ev}(f)=0$. Thus, if we take $g\in \mc D(\G)$ together with its corresponding charge $\Lambda(g)$, we get that
	$$
	0=\Scal{{\mathit ev}(f)}{\Lambda(g)}=\Scal{\Lambda(g)}{f}=\int_\G fg \,dx\,\quad\forall\, g\in \mc D(\G).
	$$
	Since $f\in L_{\mathrm loc}^{1}(\G)$, it follows that $f(x)=0$ for $\mathscr L^n$-a.e.  $x\in \G$. Therefore, $f$ turns out to be identically zero (as a function in $BV^{{{Q}/{Q-1}}}_H(\G)$).
	
	To prove that ${\mathit ev}$ is surjective, we select $\alpha\in {\bf Ch}^\ast_0(\G)$. By using Proposition \ref{lambda operator} it follows that the composition $\alpha\circ \Lambda$ belongs to the space  $(L^Q(\G))^\ast$. Hence, by the Riesz representation theorem there exists a unique $h\in L ^{{{Q}/{Q-1}}} (\G) $ for which $$\Scal{\alpha}{\Lambda(f)}=\Scal{\alpha\circ \Lambda}{f}=\int_{\G} h f\,dx\qquad\forall\, f\in L^Q(\G).$$
	We need to show that $h\in  BV^{{{Q}/{Q-1}}}_H(\G)$. To this aim, we apply the preceding  equality to ${\rm div}_H \Phi$,  whenever  $\Phi\in\mc D(\G,H\G)$ and $\|\Phi\|_\infty\le 1$. Thus $$\Scal{\alpha}{\Lambda({\rm div}_H \Phi)}=\int_{\G} h\, {\rm div}_H \Phi\,dx\,.$$ Hence,  we get that
	\begin{eqnarray*}
\int_{\G} h \, {\rm div}_H \Phi \,dx =\Scal{\alpha}{\Lambda({\rm div}_H\Phi)}\leq \|\alpha\|_{{\bf Ch}^\ast_0}\|\Lambda({\rm div}_H\Phi)\|_{{\bf Ch}_0}\leq\|\alpha\|_{{\bf Ch}^\ast_0},
	\end{eqnarray*} where the last inequality follows from \eqref{norma_lambda}.
	Taking the supremum on the left hand side  over all $\Phi\in\mc D(\G,H\G)$ such that $\|\Phi\|_\infty\le 1$, we get that $\|D_H h\|<+\infty$. If follows that $h\in  BV^{{{Q}/{Q-1}}}_H(\G)$ and that  $$\Scal{{\mathit ev}(h)}{\Lambda(f)} = \Scal{\Lambda(f)}{h}=\int_{\G} h f\,dx=\Scal{\alpha}{\Lambda(f)},$$for every $f\in L^Q(\G)$. Using that $\mc R(\Lambda)$ is dense in ${\bf Ch}_0(\G)$, we finally get that ${\mathit ev}(h)=\alpha$, as wished.\end{proof}
Notice that the map $ ev$ is  in fact an isomorphism of Banach spaces. 

\section{Bourgain-Brezis's duality argument for the getting the estimate \eqref{moonens-piconintro}}

In order to prove  inequality \eqref{1in}, Bourgain and Brezis  pass from an operator to its adjoint and conversely. A similar method is used in \cite{DePauw-Pfeffer}, \cite{DePauw-Torres}, \cite{Moonens-Picon}, and  \cite{Moonens-Picon2}.

To begin with, if  $f\in L^Q(\G)$  we have to explain in which sense we want to solve the equation  
$$
\mathrm{div}_H\Phi=f
$$
in our setting, finding a solution $\Phi\in C_0(\G, H\G)$ such that 
$$  \|\Phi\|_{L^\infty}\leq \mathscr C(Q) \|f\|_{L^Q},$$where $\mathscr C(Q)$ is a geometric constant.

The  results in this section generalize both Theorem 6.1 in \cite{DePauw-Torres} and Theorem 3.1 in \cite{Moonens-Picon} to sub-Riemannian Carnot groups.

	\subsection{A charge associated with a divergence operator}
	Also in  Carnot groups, we can  define the notion of {\it flux}. More precisely, we say that a distribution $F\in\mc D'(\G)$ is  
a flux if
the equation $\mathrm{div}_H\,Y=F$ has a continuous solution, i.e., if there exists a horizontal vector field
$Y\in C(\G;H\G)$ such that 
$$F(\varphi)=-\int_{\G}
\langle Y(x) , D_H\varphi(x)\rangle\,dx\qquad \forall  \,\varphi\in\mc D(\G).$$ 

We now have to define a linear operator $\Gamma: C_0(\G, H\G)\longrightarrow {\bf Ch}_0(\G)$ such that the charge $ \Gamma(\Phi)$,  for any given $\Phi\in C_0(\G, H\G)$, can be thought of as the (distributional) horizontal divergence of $\Phi$.

 We start by observing  that for any $f\in BV^{{{Q}/{Q-1}}}_H(\G)\subset BV_{H, {loc}}(\G)$  the structure theorem implies that 
\begin{equation}\label{riesz2-2}
\int_\G f\,\mathrm{div}_H \Phi \  dx=-\int_\G \scal{\Phi}{d[D_Hf]} \qquad\forall\, \Phi\in \mc D(\G, H\G).
\end{equation}

	We give the following definition.
	
	\begin{definition} For any $\Phi\in C_0(\G, H\G)$,  let   $\Gamma(\Phi):BV^{{{Q}/{Q-1}}}_H(\G)\longrightarrow\R$ be the linear functional defined as $$\qquad\qquad\qquad \scal{\Gamma(\Phi)}{g}:=-\int_{\G}\langle\Phi, d[D_H g]\rangle\qquad \forall\, g\in BV^{{{Q}/{Q-1}}}_H(\G).$$
	\end{definition}

	\begin{proposition}\label{gamma prop} If $\Phi\in C_0(\G, H\G)$, then $\Gamma(\Phi)\in {\bf Ch}_0(\G)$ and  $\|\Gamma(\Phi)\|_{{\bf Ch}_0}\leq  \|\Phi\|_{\infty}$.

	 As a consequence, the   linear operator  $$\Gamma:   C_0(\G, H\G)\longrightarrow {\bf Ch}_0(\G)$$ is a bounded linear operator. 
	\end{proposition}
	\begin{proof}Let $\Phi\in C_0(\G, H\G)$ and let $\{g_j\}_{j\in\N}\subset BV^{{{Q}/{Q-1}}}_H(\G)$ be a sequence such that $g_j \twoheadrightarrow 0$ as   $j\to +\infty$.  For any $\epsilon >0$,  let  $\Psi\in \mc D(\G, H\G)$ be   such that $\|\Phi-\Psi\|_{\infty}<\epsilon$. Moreover, let us set  $\mathscr K:=\sup_j \|D_H g_j\|$.  We have
		\begin{eqnarray*}\left|\Scal{\Gamma(\Phi)}{g_j}\right|&\leq&  \left|\int_{\G}\langle(\Phi-\Psi), d[D_H g_j]\rangle\,dx\right|+ \left|\int_{\G}{\rm div}_H \Psi\, g_j\,dx\right|\\&\leq&   \mathscr K \epsilon+ \left|\int_{\G}{\rm div}_H \Psi\, g_j\,dx\right|.
		\end{eqnarray*}Since ${\rm div}_H\Psi$ is a smooth compactly supported function, we get that  ${\rm div}_H\Psi\in L^Q(\G)$ and hence the second integral goes to $0$ as   $j\to+\infty$. As a consequence  $$\limsup_{j\to+\infty}\left|\Scal{\Gamma(\Phi)}{g_j}\right|\leq   \mathscr K \epsilon.$$Thus, the first claim follows from the arbitrariness of $\epsilon>0$.
		
		It is also clear that for any $g\in   BV^{{{Q}/{Q-1}}}_H(\G)$  the following inequality holds
		
		$$\left|\Scal{\Gamma(\Phi)}{g}\right| \leq \|\Phi\|_{\infty}\,\|D_H g\|.$$This    implies the second claim and achieves the proof.
	 
	\end{proof}
	
	\begin{remark}\label{=} For any $\Phi\in \mc D(\G, H\G)$, let us consider the charges $\Lambda({\rm div}_H\Phi)$ and $\Gamma(\Phi)$. It is immediate to see that $$\Scal{\Lambda({\rm div}_H\Phi)}{g}=\int_{\G}g \,{\rm div}_H\Phi \,dx= \Scal{\Gamma(\Phi)}{g}\qquad\forall\,g\in  \mc D(\G).$$Thus, using  {Theorem} \ref{density} we get that they coincide as functionals on $BV^{{{Q}/{Q-1}}}_H(\G)$.
	
Keeping this in mind, what we shall  prove in Corollary \ref{cor main} is that for any $f\in L^Q(\G)$  there exists a continuous vector field vanishing
at infinity $\Phi\in C_0(\G,H\G)$   such that
$$
\Gamma(\Phi)=\Lambda(f)
$$
in the sense that
\begin{equation}\label{significato}
	-\int_{\G}\langle\Phi, d[D_H g]\rangle=\int_{\G}f\,g\,dx \quad\quad \forall\, g\in BV^{{{Q}/{Q-1}}}_H(\G).
\end{equation}

This will be a consequence of Theorem \ref{main result} below.
\end{remark}	
Following the original idea of Bourgain and Brezis, as in \cite{DePauw-Torres} we need to characterize  the adjoint $\Gamma^\ast$ of $\Gamma$.

We first consider the map 
	$$
	-D_H: BV^{{{Q}/{Q-1}}}_H(\G)\longrightarrow \mc M(\G,H\G),
	$$
	where $D_H g:=[D_H g]$. Moreover, let $\rho:  \mc M(\G,H\G)\longrightarrow C_0(\G, H\G)^\ast$ be such that  $$\rho(\mu)(v)=\int_\G \langle v,d\mu\rangle =T_\mu(v)\qquad \forall \, v\in C_0(\G, H\G).$$
	The map $\Gamma^\ast: {\bf Ch}^\ast_0(\G)\rightarrow C_0(\G, H\G)^\ast $
makes the following diagram commutative
$$
\begin{CD}
BV^{{{Q}/{Q-1}}}_H(\G) @>-D_H>> \mathcal M(\G,H\G)\\
@V{\textit ev{}}VV     @VV {\rho}V\\
{\bf Ch}^\ast_0(\G) @>>{\Gamma^\ast}> C_0(\G, H\G)^\ast .
\end{CD}
$$
Indeed, let $\alpha\in {\bf Ch}^\ast_0(\G)$ and $\Phi\in C_0(\G, H\G)$. Furthermore, let $g={\textit ev}^{-1}(\alpha)$. Hence
\begin{eqnarray}\nonumber
	\Scal{\Gamma^\ast(\alpha)}{\Phi}_{C_0^\ast, C_0}&=&\Scal{\alpha}{\Gamma(\Phi)}_{{\bf Ch}_0^\ast, {\bf Ch}_0}
	 =\!\Scal{ev(g)}{\Gamma(\Phi)}\\\label{tutte e due}\\\nonumber&=&\!\Scal{\Gamma(\Phi)}{g}=-\int_\G  \scal{\Phi}{d[D_H g]}.
\end{eqnarray}
Thus, up to the identifications $C_0(\G,H\G)\cong \mc M(\G,H\G)$ and ${\bf Ch}^\ast_0(\G)\cong BV^{{{Q}/{Q-1}}}_H(\G)$,  since $\Gamma$ is the  distributional horizontal divergence of $\Phi$, then $\Gamma^\ast$ is (minus) the distributional  horizontal gradient $-D_Hg$ of $g$.

\begin{proposition}\label{chiusura}
The range $\mc R(\Gamma^\ast)$ of the adjoint operator  
$$
\Gamma^\ast: {\bf Ch}^\ast_0(\G) \rightarrow  C_0(\G, H\G)^\ast 
$$ 
is closed in ${\bf Ch}^\ast_0(\G).$
\end{proposition}
\begin{proof}
Let $\{\alpha_j\}_{j\in\N}\subset{\bf Ch}^\ast_0(\G)$ be a sequence such that $$\Gamma^\ast(\alpha_j)\xrightarrow[j\to+\infty]{} T,$$ for some $T\in C_0(\G,H\G)^\ast$. 
Let $\{g_j\}_{j\in\N}$ be the corresponding sequence in $BV^{{{Q}/{Q-1}}}_H(\G)$, where we have set $g_j:={\textit ev}^{-1}(\alpha_j)$ for any $j\in\N$.  The sequence $\{\Gamma^\ast(\alpha_j)\}_{j\in\N}$ is bounded, being convergent.
% {\it We claim that there exists $g\in BV^{{{Q}/{Q-1}}}_H(\G)$ such that $\mu=[D_H g]$.}
 Hence, by \eqref{tutte e due} also $\|D_Hg_j\|$  is bounded and we get that $\sup_j\|D_H f_j\|<+\infty$.  From  Proposition \ref{compact} we get that there exist $g\in BV^{{{Q}/{Q-1}}}_H(\G)$ and a  subsequence $\{g_{j_k}\}_{k\in\N}\subset BV^{{{Q}/{Q-1}}}_H(\G)$ such that $g_{j_k}\twoheadrightarrow g$ as $k\to+\infty$.	
Setting $\alpha={\textit ev}(g)$, we have 
	\begin{eqnarray*}\Scal{T}{\Phi}_{C_0^\ast, C_0}&=&\lim_{k\to+\infty}\Scal{\Gamma^\ast(\alpha_{j_k})}{\Phi}=-\lim_{k\to+\infty}\int_\G\langle\Phi, d[D_H g_{j_k}]\rangle\\
	&=& \lim_{k\to+\infty}\int_\G g_{j_k}{\rm div}_H\Phi\,dx=\int_\G g\,{\rm div}_H\Phi\,dx\\
	&=&-\int_\G
	\langle\Phi, d[D_H g]\rangle=\Scal{\Gamma^\ast(\alpha)}{\Phi}_{C_0^\ast, C_0}
	\end{eqnarray*}for any $\Phi\in \mc D(\G,H\G)$.

	 By the density of $\mc D(\G,H\G)$  
 in $C_0(\G,H\G)$, we get that $T=\Gamma^\ast(\alpha)$, which achieves the proof.

\end{proof}

As a corollary, keeping in mind Proposition II.18 in \cite{{Brezis}} we have the following:
\begin{corollary}\label{cor range}
The range $\mc R(\Gamma)$ of $\Gamma$ is closed in ${\bf Ch}_0(\G)$.
\end{corollary}
\subsection{Main results}
We are in a position to solve the problem 
$$
\Gamma(\Phi)=F,
$$
whenever $\Phi\in C_0(\G, H\G)$ and $F\in\mc D'(\G)$.  More precisely,  the following holds:
%\begin{proof}
%
%\end{proof}
%Thanks to the following theorem we are going to obtain that if $F\in\mc D'(\G)$, then, there exists $\Phi\in C_0(\G,H\G)$ such that \begin{equation}\label{eq div}\Gamma(\Phi)=F\end{equation} if, and only if,  $F\in{\bf Ch}_0(\G)$. The necessity part is given by Proposition \ref{gamma prop}. 
\begin{theorem}\label{main result}
Let $F\in\mc D'(\G)$. Then, there exists $\Phi\in C_0(\G,H\G)$ such that \begin{equation}\label{eq div}\Gamma(\Phi)=F\end{equation} if, and only if,  $F\in{\bf Ch}_0(\G)$.  

In addition, if $F\in{\bf Ch}_0(\G)$ there exists a solution $\Phi\in C_0(\G,H\G)$ of \eqref{eq div} such that
\begin{equation}
\label{pre-moonens-picon} \|\Phi\|_\infty\leq 2  \|F\|_{{\bf Ch}_0}\,.
\end{equation}
\end{theorem}

\begin{proof} 
%The proof is divided into two main steps, and each step in a number of clai	

	\noindent{\bf Step 1 (proof of \eqref{eq div}).\,}  The necessity part follows from Proposition \ref{gamma prop}. 
	
	Furthermore, since in Corollary \ref{cor range} we have proved that $\mc R(\Gamma)$ is closed, the sufficiency part will be proved once we have shown  that $\mc R(\Gamma)$ is dense in ${\bf Ch}_0(\G)$.
	
	To show  that  $\mc R(\Gamma)$ is dense in ${\bf Ch}_0(\G)$ we  use a standard
	consequence of the Hahn-Banach theorem; see \cite{Brezis}, Corollary I.8. We assume that $\alpha\in {\bf Ch}^\ast_0(\G)$  vanishes 
	on all of $\mc R(\Gamma)$. Thus, we have to show that  $\alpha$ must vanish everywhere on ${\bf Ch}_0(\G)$. 
	
To this aim, let  $\alpha\in {\bf Ch}^\ast_0(\G)$ be such that $\Scal{\alpha}{\Gamma(\Phi)}=0$ for every $\Phi\in C_0(\G,H\G)$. By Proposition \ref{isom},  there exists a unique $g\in BV^{{{Q}/{Q-1}}}_H(\G)$ such that $\alpha={\mathit ev}(g)$. Then $$0= \Scal{{\mathit ev}(g)}{\Gamma(\Phi)}=  \Scal{\Gamma(\Phi)}{g}= -\int_{\G}\langle \Phi, d[D_H g]\rangle\qquad\forall\,\Phi\in C_0(\G,H\G).$$ 
	This implies that $D_Hg=0$ and in turn that $g=0$, since $g\in BV^{{{Q}/{Q-1}}}_H(\G)$.
	%(in fact, $f$ must be constant on the whole of $\G$ because its variation measure vanishes, but a constant function belongs to $L^{Q/Q-1}(\G)$ only if it is a.e. $0$). \\

	\noindent{\bf Step 2 (proof of the second part).\,} Let $F\in {\bf Ch}_0(\G)$. We show that  it is possible to find a solution of \eqref{eq div} that satisfies {\it also} the estimate \eqref{pre-moonens-picon}.

	Again, we use the original idea by Bourgain and Brezis for periodic functions, already  used in the Euclidean setting by \cite{DePauw-Pfeffer} and \cite{{Moonens-Picon}}, under more general assumptions.
	
	For the sake of simplicity,  we will set  here $X= C_0(\G,H\G)$. 
	
	Let $F\in {\bf Ch}_0(\G)$ be such that $\|F\|_{{\bf Ch}_0}>0$, and define two convex subsets by setting $${\mc U}:=\left\{ \Phi\in X\;:\; \Gamma(\Phi)=F\right\},\qquad {\mc V}:=\left\{ \Phi\in X\;:\; \|\Phi\|_{\infty}< 2\|F\|_{{\bf Ch}_0}\right\}.$$
	%where ${\mathscr C}_{_{\! GN}}$ is the Gagliardo-Nirenberg constant.  
From Step 1 we get that ${\mc U}\neq \emptyset$. Moreover, ${\mc V}\neq \emptyset$ because $\Phi=0$  clearly  belongs to ${\mc V}$. 
	
	\noindent{\bf Claim:\,} {\it We claim that ${\mc U}\cap {\mc V}\neq \emptyset$.} 
	
	If we could show that the claim is true, then 
		 the proof would be complete since we would have found a solution $\Phi$ of \eqref{eq div} that satisfies  also the estimate $\|\Phi\|_{\infty}< 2\|F\|_{{\bf Ch}_0}$. 
		%Therefore, using Proposition \ref{lambda operator} $\|\Phi\|_{\infty}< 2\|\Lambda(f)\|_{{\bf Ch}_0}\le 2{\mathscr C}_{_{\! GN}}\|f\|_{L^Q}$ we get eventually \eqref{moonens-picon} (remember, ${\mathscr C}_{_{\! GN}}$ is the Gagliardo-Nirenberg constant appearing in \eqref{gagliardo}).\\
		
	Thus,	we are  left to prove the claim.
	By contradiction, we assume that \begin{equation}\label{assurdo} {\mc U}\cap {\mc V}=\emptyset.
	\end{equation}By the first geometric form of the Hahn-Banach theorem (see, e.g.,  Theorem 1.6 in \cite{Brezis}) we get that there exist $T\in X^\ast$ and $t\in\R$ such that:
	\begin{equation}\label{H-B}
		\Scal{T}{\Phi}\geq t\quad\forall\, \Phi\in {\mc U}\quad\mbox{and}\quad \Scal{T}{\Phi}\leq t\quad\forall\, \Phi\in {\mc V}.
	\end{equation}
Note that $t>0$, since $\Phi=0\in \mc V$.    { Moreover, we observe that ${\rm Ker}(\Gamma)\subset {\rm Ker}(T)$.}

	In fact, let $\Phi_0\in  {\rm Ker}(\Gamma) $ and   $\Phi\in \mc U$.  Then, for every $s\in\R$ we have $\Phi+s\Phi_0\in \mc U$. \\As a consequence, from the inequality $ \Scal{T}{\Phi+s\Phi_0 }\ge t$ we should have $$ s \Scal{T}{ \Phi_0 }\ge t- \Scal{T}{\Phi}\qquad\forall\,s\in\R.$$ But this   does not hold unless $\Scal{T}{ \Phi_0 }=0$. Hence $\Phi_0\in{\rm Ker}(T)$.  
	 Being surjective, $\Gamma$ is also open by the open mapping theorem. Therefore, it turns out that $\Gamma$ is a quotient map. Hence there exists $\alpha\in {\bf Ch}_0(\G)^\ast$ such that $T=\alpha\circ \Gamma$. Now, take $\tilde g={\textit ev}^{-1}(\alpha)\in BV^{{{Q}/{Q-1}}}_H(\G)$. 
	Then,  for any $\Phi\in X$ we have \begin{eqnarray}\nonumber
		-\int_\G  \scal{\Phi}{d[D_H \tilde g]}&=&\Scal{\Gamma(\Phi)}{\tilde g}=\Scal{{\textit ev}(\tilde g)}{\Gamma(\Phi)}=\Scal{\alpha}{\Gamma(\Phi)}_{{\bf Ch}_0^\ast,{\bf Ch}_0}\\ 	\label{espressione T}\\\nonumber &=&(\alpha\circ \Gamma)(\Phi)=\Scal{T}{\Phi}_{C_0^\ast,C_0}.
	\end{eqnarray}

	On the other hand, let $\Phi\in{\mc D}(\G, H\G)$ be such that  $\|\Phi\|_{\infty}\le 1$ and  choose $\epsilon >0$ such that $1+\varepsilon <2$. Hence $\Psi:=({1+\varepsilon})\|F\|_{{\bf Ch}_0}{\Phi}\in \mc V$. In addition,  we have
	\begin{eqnarray*}
		\int_{\G} \tilde g\, {\rm div}_H\Phi\,dx&=&-\int_{\G}\langle\Phi, d[D_H\tilde g]\rangle=-\frac{1}{({1+\varepsilon})\|F\|_{{\bf Ch}_0}}\int_{\G}\langle\Psi, d[D_H\tilde g]\rangle\\
		&\underset{\rm by \ \eqref{espressione  T}}{=}&\frac{1}{({1+\varepsilon})\|F\|_{{\bf Ch}_0}}\Scal{T}{\Psi}\underset{\Psi\in\mc V}{\le} \frac{t}{({1+\varepsilon})\|F\|_{{\bf Ch}_0}}\,.
	\end{eqnarray*}
	In particular, by taking  the supremum on all $\Phi\in{\mc D}(\G, H\G)\subset X$ such that $\|\Phi\|_{\infty}\le 1$, we get that
	$$
	\|D_H\tilde g\|\le \frac{t}{({1+\varepsilon})\|F\|_{{\bf Ch}_0}}.
	$$
Let $\Phi\in \mc U$. 	Using the last estimate together with \eqref{H-B} and \eqref{espressione T}, we get that
	\begin{eqnarray*}
		t\leq \Scal{T}{\Phi}= \Scal{\Gamma(\Phi)}{\tilde g} \underset{\Phi\in \mc U}{=}\Scal{F}{\tilde g}\le \|F\|_{{\bf Ch}_0} \|D_H\tilde g\|\le \frac{t}{({1+\varepsilon})}.
	\end{eqnarray*} But this cannot be true, since we have seen that $t$ is positive.  This contradiction shows our claim and concludes the proof.

	\end{proof}

By Step 1 of the above proof, we have that $\Gamma $ is surjective. Since $\Gamma$ is  also continuous (see Proposition \ref{gamma prop}), by the open mapping theorem there exists a positive constant $\mathscr C>0$ such  that 
$$
 \|\Phi\|_\infty\leq {\mathscr C}  \|F\|_{{\bf Ch}_0},
$$for any  solution $\Phi\in C_0(\G,H\G)$ of \eqref{eq div}. 
Therefore, the second part of  Theorem \ref{main result} would follow straightforwardly {\it for any solution} $\Phi\in C_0(\G,H\G)$  of \eqref{eq div}, if one were satisfied with a generic constant. On the contrary,  we have been able to get an estimate  with an explicit constant,  but paying the price that the estimate holds  for {\it some} $\Phi$. We  also note that the constant $2$ does not play any role here. The proof would work  as well  with a constant as close to $1$ as one wants.  
\\

As an immediate corollary of Theorem \ref{main result}, for any $f\in L^Q(\G)$ we have the following estimate  with a geometric constant, which depends only on the homogeneous dimension. 

Clearly, the equation $\Gamma(\Phi)=\Lambda(f)$ is meant here as specified in \eqref{significato}.
\begin{corollary}\label{cor main}
For any $f\in L^Q(\G)$  there exists a  solution $\Phi\in C_0(\G,H\G)$ of 
\begin{equation}\label{eq div2}\Gamma(\Phi)=\Lambda(f)\end{equation}
satisfying the  inequality  
\begin{equation}
\label{moonens-picon} \|\Phi\|_\infty\leq 2 \ {{\mathscr C}_{_{\! GN}}} \|f\|_{L^Q},
\end{equation}where ${\mathscr C}_{_{\! GN}}$ is the constant appearing in \eqref{gagliardo-bis}.
\end{corollary}
\begin{proof}
Recall that  $\Lambda(f)\in {\bf Ch}_0(\G)$ for any $f\in L^Q(\G)$. Thus, from Theorem \ref{main result} we get that there exists a solution  $\Phi\in C_0(\G,H\G)$ satisfying
$$\|\Phi\|_\infty\leq 2  \|\Lambda(f)\|_{{\bf Ch}_0}.$$
Finally,  \eqref{moonens-picon} follows from Proposition \ref{lambda operator}.
\end{proof}
 
	%
	%
	%
	%
	%
		%\noindent{\bf Claim 2.\,} {\it We claim that  $\mc R(\Gamma)$ is closed in ${\bf Ch}_0(\G)$.}
		%
	%Observe that  by the Closed Range Theorem we only need to show that the range $\mc R(\Gamma^\ast)$  of the adjoint operator $\Gamma^\ast$ is closed in $C_0(\G,H\G)$.  
	%
	%
	%
%In order to prove Claim 2, we preliminary characterize the adjoint operator $\Gamma^\ast$. 
	%
	%
		%\noindent{\bf Claim 2.1.\,} {\it We claim that  $\Gamma^\ast\circ{\mathit ev}=-D_H$.} 
		%
		%To see this, we  use  the two identifications $$C_0(\G,H\G)\cong \mc M(\G,H\G),\qquad {\bf Ch}^\ast_0(\G)\cong BV^{{{Q}/{Q-1}}}_H(\G).$$Thus, for any $\alpha\in  {\bf Ch}^\ast_0(\G)$ there exists  one, and only one, $f\in BV^{{{Q}/{Q-1}}}_H(\G)$ such that $\alpha={\mathit ev}_f$. As a consequence  $$\Scal{\Gamma^\ast({\mathit ev}_f)}{\Phi}= \Scal{{\mathit ev}_f}{\Gamma(\Phi)}= \Scal{\Gamma(\Phi)}{f}=-\int_{\G}\langle\Phi, d[D_Hf]\rangle,$$which implies that $$\Gamma^\ast\circ {\mathit ev}=-D_H,$$as claimed. 
		%
		%
		%
	%
	%
	%
%Now, 
%

\section*{Funding informations}
{\small

A.B. is supported by the University of Bologna, funds for selected research topics, and by MAnET Marie Curie
Initial Training Network, and by 
GNAMPA of INdAM (Istituto Nazionale di Alta Matematica ``F. Severi''), Italy.\\F.M. is supported by 
GNAMPA of INdAM (Istituto Nazionale di Alta Matematica ``F. Severi''), Italy.
}

\bibliographystyle{amsplain}

\vspace{0.5cm}
\bigskip
\tiny{
\noindent
Annalisa Baldi: 
\par\noindent
Universit\`a di Bologna,\\ Dipartimento
di Matematica\par\noindent Piazza di
Porta S.~Donato 5, 40126 Bologna, Italy 
\par\noindent
e-mail:
annalisa.baldi2@unibo.it 
}

\bigskip

\tiny{
\noindent
Francescopaolo Montefalcone:
\par\noindent  Universit\`a di Padova,\\ Dipartimento
di Matematica ``Tullio Levi-Civita'' \par\noindent Via Trieste 63, 35121 Padova,  Italy 
\par\noindent
 e-mail: montefal@math.unipd.it
}

\end{document}